\documentclass[10pt,leqno]{amsart}
\usepackage{graphicx} 

\DeclareUnicodeCharacter{2161}{\^{}}

\usepackage{eucal}

\usepackage[backend=biber,style=alphabetic,doi=true,isbn=true, maxalphanames=4, 
  minalphanames=4,maxbibnames=99, url = false]{biblatex}
\usepackage[english]{babel}
\usepackage{microtype}
\usepackage{mathrsfs}
\usepackage{amsthm, thmtools}
\usepackage{amsmath, centernot}
\usepackage{amssymb}
\usepackage{indentfirst,csquotes}
\usepackage{amscd}
\usepackage{mathtools}
\usepackage{wasysym}
\usepackage{braket}
\usepackage{appendix}
\usepackage{dsfont}
\usepackage{caption}
\usepackage{subfig}
\usepackage[all]{xy}
\usepackage{xfrac}
\usepackage{hyperref}
\usepackage{color}
\usepackage{bbm}
\usepackage{stmaryrd}
\usepackage{fancyhdr}
\usepackage{graphicx}
\usepackage{nomencl}
\usepackage{ esint }
\usepackage{enumitem}

\numberwithin{equation}{section}

\allowdisplaybreaks

\newtheorem{teorema}{Theorem}[section]
\newtheorem{prop}[teorema]{Proposition}
\newtheorem{co}[teorema]{Corollary}
\newtheorem{lemma}[teorema]{Lemma}

\newtheorem{df}[teorema]{Definition}

\newtheorem{oss}[teorema]{Remark}

\newtheorem{ass}[teorema]{Assumption}

\newcommand{\R}{\mathbb{R}}

\newcommand{\Q}{\mathbb{Q}}
\newcommand{\N}{\mathbb{N}}

\newcommand{\PP}{\mathcal{P}}

\newcommand{\X}{\textit{X}}

\newcommand{\e}{\mathrm{e}}

\newcommand{\supp}{\operatorname{supp}}

\newcommand{\KFP}{\operatorname{KFP}}
\newcommand{\MP}{\operatorname{MP}}

\newcommand{\nablasquare}{\nabla^{\otimes_2}}

\title[First order equation on random measures]{First order equation on random measures as superposition of weak solutions to the McKean-Vlasov equation}
\author{
  Alessandro Pinzi
 }\thanks{Bocconi University, Department of Decision Sciences, via Roentgen 1, 20136 Milano (Italy). 
Email: \textsf{alessandro.pinzi@phd.unibocconi.it}.  \textit{Orcid}: \textsf{https://orcid.org/0009-0007-9146-5434}}
\date{}

\begin{document}

\begin{abstract}
The goal of this paper is to define an evolution equation for a curve of random probability measures $(M_t)_{t\in[0,T]}\subset \PP(\PP(\R^d))$ associated to a non-local drift $b:[0,T]\times\R^d \times \PP(\R^d) \to \R^d$ and a non-local diffusion term $a:[0,T]\times \R^d \times \PP(\R^d) \to \operatorname{Sym}_+(\R^{d\times d})$. Then, we show that any solution to such an equation on random measures can be lifted twice: to a superposition of solutions to a non-linear Kolmogorov-Fokker-Planck equation and to a superposition of weak solutions to the McKean-Vlasov equations. Finally, we exploit this nested superposition result to show how existence and uniqueness can be transferred from the equation on random measures to the associated non-linear Kolmogorov-Fokker-Planck equation and to the McKean-Vlasov equation, assuming uniqueness of the linearized version of KFP. As a tool, we will introduce integral metrics over the spaces of probability measures $\PP(\R^d)$ in duality with smooth functions, including a weighted second-order metric.
\end{abstract}

\maketitle

{\small
		\par
	}


\section{Introduction}
The superposition principle plays an important role in many evolution problems: its first version was proved by L. Ambrosio (see \cite[Theorem 8.2.1]{ambrosio2005gradient}), and relates the Eulerian (continuity equation) and the Lagrangian description (system of ODEs) of the flow led by a vector field $v:[0,T]\times \R^d \to \R^d$. It was introduced to study the well-posedness of ODEs (in a selection sense) under non-smooth assumptions on the vector field \cite{ambrosio2004transport, ambcri06}, extending the celebrated work by R. Di Perna and P.L. Lions \cite{diperna1989cauchy}. It is a useful tool in many other contexts (e.g. optimal transport), and for this reason, it has been extended to more abstract spaces (see e.g. \cite{stepanov2017three}).

Let us recall the following two versions of the superposition principle:
\begin{itemize}
    \item[(A)] the \textit{stochastic superposition principle} proved in \cite{figalli2008existence, trevisan2016well, bogachev2021ambrosio}, that in the spirit of the Ambrosio's superposition principle, relates weak solutions of a stochastic equation to the associated linear Kolmogorov-Fokker-Planck equation (see also §\ref{sec: fin dim equations});
    \item[(B)] the \textit{nested superposition principle} proved in \cite{pinzisavare2025}, where the authors defined an abstract continuity equation for a curve of random measures $(M_t)_{t\in[0,T]}\subset \PP(\PP(\R^d))$ that can be written either as superposition of solutions to a non-local continuity equation $\partial_t\mu_t + \operatorname{div}(b_t(\cdot,\mu_t)\mu_t) = 0 $ or as superposition of solutions to an interacting particle systems $dX_t = b_t(X_t,\operatorname{Law}(X_t))dt$.
\end{itemize}

In this paper, we want to prove a new nested superposition principle for random measures (which in particular is an extension of (B)) that relies on the stochastic superposition principle (A). To be more specific, consider a non-local vector field $b:[0,T]\times \R^d \times \PP(\R^d) \to \R^d$ and a non-local diffusion term $\sigma:[0,T]\times \R^d \times \PP(\R^d) \to \R^{d\times m}$, and define $a:= \sigma \sigma^\top$. The term non-local highlights their dependence on the measure variable. For any $\mu\in \PP(\R^d)$, define the operator
\[L^{\mu}_{b_t,a_t}\phi(x) := \sum_{i=1}^d b_t^{i}(x,\mu)\cdot \partial_i\phi(x) +\frac{1}{2}\sum_{i,j=1}^d a^{i,j}_t(x,\mu)\partial_{i,j}\phi(x), \]
for all $\phi\in C_b^2(\R^d)$.
The stochastic superposition principle shows that there is a correspondence (possibly non 1-1) between curves of probability measures $(\mu_t)_{t\in[0,T]} \subset \PP(\R^d)$ that solve (in a distributional sense) the non-linear Kolmogorov-Fokker-Planck equation
\begin{equation}\label{eq: intro KFP}
    \partial_t\mu_t = (L_{b_t,a_t}^{\mu_t})^*\mu_t,
\end{equation}
and weak solutions to the McKean-Vlasov equation (or equivalently solutions of the associated non-linear martingale problem \cite{stroock1997multidimensional})
\begin{equation}\label{eq: intro Mckean vlasov}
    dX_t = b_t(X_t,\operatorname{Law}(X_t))dt + \sigma_t(X_t,\operatorname{Law}(X_t))dW_t.
\end{equation}

We will introduce an evolution equation over random measures $\PP(\PP(\R^d))$ that can be seen either as a superposition of solutions to \eqref{eq: intro KFP} or as a superposition of weak solutions to \eqref{eq: intro Mckean vlasov}. To this aim, we define an operator $\mathcal{K}_{b_t,a_t}$ acting on cylinder functions $\operatorname{Cyl}_b^{1,2}(\PP(\R^d))$ (see Definition \ref{def: cyl functions}, they will play the role of smooth test functions), so that the equation $\partial_tM_t = \mathcal{K}_{b_t,a_t}^*M_t$ is well defined (see Definition \ref{def: sol of nd KFP for rm}), where $(M_t)_{t\in[0,T]}\subset \PP(\PP(\R^d))$. Then, the main result of this paper can be summarized as follows.

\begin{teorema}\label{main theorem}
    Let $b:[0,T]\times \R^d \times \PP(\R^d) \to \R^d$ and $a:[0,T]\times \R^d \times \PP(\R^d) \to \operatorname{Sym}_+(\R^{d\times d})$ be Borel measurable maps. Let $\boldsymbol{M} = (M_t)_{t\in[0,T]} \in C([0,T],\PP(\PP(\R^d)))$ be such that \begin{equation}\label{eq: equation main theorem}
    \partial_tM_t = \mathcal{K}_{b_t,a_t}^*M_t \quad \text{and}\quad \int_0^T \int_{\PP(\R^d)}\int_{\R^d} \frac{|b_t(x,\mu)|}{1+|x|} + \frac{|a_t(x,\mu)|}{1+|x|^2} d\mu(x)dM_t(\mu)dt <+\infty.
    \end{equation}
    Then, there exist $\Lambda \in 
    \PP(C([0,T],\PP(\R^d)))$ and $\mathfrak{L}\in \PP(\PP(C([0,T],\R^d)))$ liftings of $\boldsymbol{M}$, in the sense that $(\mathfrak{e}_t)_\sharp \Lambda = M_t$ and $(E_t)_\sharp \mathfrak{L} = M_t$ (see \eqref{eq: maps for marginals}) for all $t\in[0,T]$ satisfying:
    \begin{enumerate}
        \item $\Lambda$-a.e. $\boldsymbol{\mu} = (\mu_t)_{t\in[0,T]} \in C([0,T],\PP(\R^d))$ solves \eqref{eq: intro KFP} and it satisfies 
        \begin{equation}\label{eq: integr of Lambda}
        \int \int_0^T \int_{\R^d} \frac{|b_t(x,\mu_t)|}{1+|x|} + \frac{|a_t(x,\mu_t)|}{1+|x|^2} d\mu_t(x)dt d\Lambda(\boldsymbol{\mu}) <+\infty;
        \end{equation}
        \item $\mathfrak{L}$-a.e. $\lambda \in \PP(C_T(\R^d))$ is a martingale solution of \eqref{eq: intro Mckean vlasov} (see Definition \ref{def: MP solution fin dim}) and it satisfies
        \begin{equation}\label{eq: integr of mathfrak L}
        \int  \int \int_0^T\frac{|b_t(\gamma_t,(\e_t)_\sharp \lambda)|}{1+|\gamma_t|} + \frac{|a_t(\gamma_t,(\e_t)_\sharp \lambda)|}{1+|\gamma_t|^2} dt d\lambda(\gamma) d\mathfrak{L}(\lambda) <+\infty;
        \end{equation}
        \item $E_\sharp \mathfrak{L} = \Lambda$ (see \eqref{eq: maps for marginals}) and there exists a map $G_{b,a}: C([0,T],\PP(\R^d)) \to \PP(C([0,T],\R^d))$ such that $(G_{b,a})_\sharp \Lambda = \mathfrak{L}$ and $E(G_{b,a}(\boldsymbol{\mu})) = \boldsymbol{\mu}$ for 
        $\Lambda$-a.e. $\boldsymbol{\mu}$.
    \end{enumerate}
    Conversely: 
    \begin{enumerate}
        \item[(i)] given $\mathfrak{L}\in \PP(\PP(C([0,T],\R^d)))$ satisfying the condition in (2), then $E_\sharp \mathfrak{L} \in \PP(C([0,T],\PP(\R^d)))$ satisfies the condition in (1);
        \item[(ii)] given $\Lambda \in \PP(C([0,T],\PP(\R^d)))$ satisfying the condition in (1), then the curve of random measures defined by $M_t:= (\mathfrak{e}_t)_\sharp \Lambda$, satisfies \eqref{eq: equation main theorem}. 
    \end{enumerate} 
\end{teorema}

This theorem is followed by a uniqueness scheme: indeed, under suitable assumptions, we can prove that uniqueness for the equation on random measures implies uniqueness of \eqref{eq: intro KFP} and \eqref{eq: intro Mckean vlasov} (see §\ref{sec: uniqueness}). 

As a tool to prove the main theorem, we introduce complete integral metrics over the set of non-negative and finite measures $\mathcal{M}_+(\R^d)$ (and in particular over $\PP(\R^d)$) in duality with smooth functions, which we prove to induce the narrow topology. 

The advantages of this approach are many:
\begin{itemize}
    \item the equation $\partial_t M_t = \mathcal{K}_{b_t,a_t}^*M_t$ shares interesting properties. Firstly, it is linear in $M_{\cdot}$, making it an infinite-dimensional linearized version of the Kolmogorov-Fokker-Planck equation. Secondly, it is a first order equation, in the sense that the operator $\mathcal{K}_{b_t,a_t}$ satisfies the Leibniz rule $\mathcal{K}_{b_t,a_t}(FG) = F\mathcal{K}_{b_t,a_t}G + G \mathcal{K}_{b_t,a_t} F$;
    \item in principle, it can be used to study well-posedness for \eqref{eq: intro KFP} and \eqref{eq: intro Mckean vlasov}, in a selection sense, under low-regularity for the coefficients, in analogy with the original works by Ambrosio (for the deterministic case) and Figalli-Trevisan (for the stochastic case).
\end{itemize}

Of course, there is no free lunch: the first thing to do to pursue the previous way for well-posedness is to find a good class of curves of random measures in which to look for existence and uniqueness, under suitable assumptions on the coefficients, that must be understood as well. To find such a class seems to be a challenging problem, due to the infinite-dimensionality of the space $\PP(\R^d)$ and the lack of a nice reference measure on it, which prevents the application of standard finite-dimensional techniques. This aspect is worth further investigating, for example endowing the space $\PP(\R^d)$ with a Gaussian-regular measure (see \cite{PS25convex}).

The well-posedness for McKean-Vlasov equation attracted many researchers in recent years \cite{barbu2024nonlinear,chaintron2022propagation}. In \cite{barbu2024nonlinear}, the authors grouped the results of a series of papers in which they studied the McKean-Vlasov equation, mainly exploiting the theory developed by A. Figalli and D. Trevisan. In particular, in \cite{barbu2021uniqueness}, they showed how the uniqueness of a linearized version of the McKean-Vlasov equation is useful to gather uniqueness for the original equation. In §5, Assumption \ref{assumption uniqueness} is a fundamental assumption on uniqueness of the linearized McKean-Vlasov, and it is exploited as in \cite{barbu2021uniqueness} to prove the uniqueness scheme.

Other approaches for the well-posedness of the McKean-Vlasov equation (recovering also strong uniqueness) are through uniform continuity assumptions on the coefficients \cite{de2020strong, pascucci2025mckean, de2022well}. In particular, some Lipschitz/Hölder continuity assumptions are done w.r.t. the state and the measure variable. To do so, several distances are considered in the measure variable. In §3, we introduce a new natural distance, in duality with $C^\ell$ functions, and the case with $\ell = 2$ will be an important tool for the proof of Theorem \ref{main theorem}. It would be interesting to understand if other metric-like discrepancies of order $2$ between probability measures (e.g. \cite{huesmann2019benamou, brigati2025kinetic, bolbotowski2024kantorovich}) may give a metric characterization of the equation $\partial_tM_t  = \mathcal{K}_{b_t,a_t}^*M_t$, as in \cite[Section 3 \& 4]{pinzisavare2025}.

A similar result to our main theorem was already obtained in \cite{lacker2022superposition}, where they considered a more general equation, putting a stochastic noise at the level of equation \eqref{eq: intro KFP}, raising an operator $\mathcal{K}$, at the level of random measures, of diffusive nature. The novelties of this paper are both in the approach, that relies on a measurable selection argument (see §\ref{sec: nd equations}) and in the integrability assumptions; indeed, we only require integrability of
\begin{equation}
    \frac{b_t(x,\mu)}{1+|x|} \quad \text{ and } \quad \frac{a_t(x,\mu)}{1+|x|^2},
\end{equation}
while they could perform the liftings only under $L^p$-assumptions of $a$ and $b$, with $p>1$ (see also \cite{rehmeier2023linearization}). A fundamental aspect for this approach is the Borel measurability of the set of solution to \eqref{eq: intro KFP} and the sets of martingale solutions to \eqref{eq: intro Mckean vlasov} under very general assumptions.  Moreover, we also show the uniqueness equivalences under the uniqueness assumption for the linearized Kolmogorov-Fokker-Planck equation, that allows one to recover uniqueness under uniformly (w.r.t. the variable $\mu$) Lipschitz assumption in the variable $x$ (see Lemma \ref{lemma: unif lip assumption}).


\bigskip

\textbf{Outline of the paper.} In \textbf{Section \ref{sec: prel}}, we fix the main spaces we are going to use in the paper, in particular fixing some natural topologies on them. Then, we recall the stochastic superposition principle, introducing the linear Kolmogorov-Fokker-Planck equation and the martingale problem. It will be the notion of solution that we use for stochastic differential equations, that is equivalent to weak solutions.

In \textbf{Section \ref{sec: metric}}, we introduce the distance $D_\ell$ over probability measures, in duality with $C_0^\ell$ functions with controlled norm. We will see that they are complete metrics that induce the narrow topology and will be an important tool for our main proof, but we believe they may have their own independent interest.

In \textbf{Section \ref{sec: nd equations}}, we introduce the operator $\mathcal{K}_{b_t,a_t}$ acting on cylinder functions and associated to the non-local drift term $b:[0,T]\times \R^d \times \PP(\R^d) \to \R^d$ and the non-local diffusion term $a:[0,T]\times \R^d \times \PP(\R^d) \to \operatorname{Sym}_+(\R^{d\times d})$. Then we prove the nested stochastic superposition principle, Theorem \ref{main theorem}, in several steps: first, we show Claims (i) and (ii); then, we prove the existence of $\Lambda$ satisfying Claim (1), embedding $\PP(\R^d)$ (endowed with the distance $D_{2,w}$ introduced in Section \ref{sec: metric}) in $\R^\infty$, in which we transfer the equation to use the superposition principle in $\R^\infty$ proved in \cite{ambrosio2014well}; then, we show the existence of $\mathfrak{L}$ that satisfies Claims (2) and (3), using a measurable selection argument. The non-trivial part is to prove the measurability of the set $\operatorname{KFP}(b,a)$ and $\operatorname{MP}(b,a)$ (see Definition \ref{def: KFP(b,a) and MP(b,a)}).

Finally, in \textbf{Section \ref{sec: uniqueness}}, we show how uniqueness can be transferred from the equation on random measures to the non-linear Kolmogorov-Fokker-Planck and to weak solutions of the McKean-Vlasov equation.

\section{Preliminaries}\label{sec: prel}

\subsection{Canonical topology over spaces of measures and the space of continuous curves}
Let $(X,\tau)$ be a Polish space, i.e. for which there exists a distance $d$ that induces the topology $\tau$ and makes it a complete and separable space. The collection $\mathcal{B}(X)$ denotes the Borel $\sigma$-algebra of $X$, i.e. the $\sigma$-algebra generated by the topology $\tau$.

\subsubsection{Space of continuous curves}\label{subsub: space of curve} The space of continuous curves from $[a,b]$ to $X$, denoted by $C([a,b],X)$ is naturally endowed with the \textit{compact-open topology}. 
Its topology can be metrized as well, resulting as Polish (see \cite[Theorem 2.4.3]{srivastava2008course} for separability): indeed, considering a distance $d$ over $X$ that induces its topology $\tau$, the compact open topology is metrized by the sup-distance over curves, i.e. 
\begin{equation}
    D_d(\boldsymbol{x}_1, \boldsymbol{x}_2):= \sup_{t\in [a,b]} d(\boldsymbol{x}_1(t), \boldsymbol{x}_2(t)) \quad \forall \boldsymbol{x}_1, \boldsymbol{x}_2 \in  C([a,b];X).
\end{equation}
It will be useful, in this context, to consider a bounded distance $\hat{d}$ that induces the compact-open topology. In particular, given any distance $d$ inducing $\tau$, defining the truncated one, i.e. $\hat{d} = d\wedge 1$, we have that the compact-open topology over $C([a,b],X)$ is induced by $D_{\hat{d}}$ as well. 

Fixing a distance 
$d$ that induces the topology of $X$, the set of absolutely continuous curves denoted by $AC([a,b],X)$ is defined as the collection of continuous curves $\boldsymbol{x} = (x_t)_{t\in[a,b]}$ for which there exists a function $g:(a,b) \to [0,+\infty]$ that is in $L^1(a,b)$ and satisfying 
\[d(x_s,x_t) \leq \int_{s}^t g(r)dr \quad \text{ for all } a\leq s < t \leq b.\]
For absolutely continuous curves, there always exists the metric derivative, i.e. 
\[|\dot{\boldsymbol{x}}|_d(t) :=\lim_{h\to 0} d(x_{t+h},x_{t})\]
exists for a.e. $t\in (a,b)$. Moreover $|\dot{\boldsymbol{x}}|(\cdot) \in L^1(a,b)$ and it is (pointwise) the smallest function $g$ that can be considered in the definition.

When $a=0$, the space $C([0,b],X)$ will be denoted as $C_b(X)$.

\subsubsection{Spaces of measures}
We will denote with $\PP(X)$ the space of Borel probability measures over $X$. More generally, we denote by $\mathcal{M}_+(X)$, $\mathcal{M}(X)$ and $\mathcal{M}(X;\R^n)$, respectively, the space of finite positive measures, signed measures and $\R^n$-valued measures, in both the last two cases with finite total variation, where for a given $\nu \in \mathcal{M}(X;\R^n)$, its total variation is 
\begin{equation}
    |\nu|(A) := \sup\left\{ \sum_{n=1}^{+\infty} |\nu(A_n)| \ : \ \bigcup A_n = A, \ A_i \cap A_j = \emptyset \text{ as } i\neq j \right\}.
\end{equation}
Notice that $\PP(X) \subset \mathcal{M}_+(X) \subset \mathcal{M}(X)$. In particular, we will endow $\mathcal{M}(X;\R^n)$ with the narrow topology, so that $\PP(X)$ and $\mathcal{M}_+(X)$ are endowed with the subspace topology, and actually are closed subsets. 

The narrow topology over $\mathcal{M}(X;\R^n)$ is the smallest topology for which the functional 
\[\mathcal{M}(X;\R^n) \ni \nu \mapsto \int_X \phi(x) \cdot d\nu(x)\]
is continuous for all $\phi \in C_b(X;\R^n)$, i.e. bounded and continuous functions over $X$ taking values in $\R^n$. It is important to recall that:
\begin{itemize}
    \item[(a)] the spaces $\PP(X)$ and $\mathcal{M}_+(X)$, endowed with the narrow topology, are Polish;
    \item[(b)] the space $\mathcal{M}(X;\R^n)$ is not Polish, but is still a Lusin space (see e.g. \cite[Remark 2.4]{pinzisavare2025}).
\end{itemize}

Recall the push-forward operation defined over positive measures: if $f:(Z_1,\mathcal{F}_1) \to (Z_2,\mathcal{F}_2)$ is a function between two generic measurable spaces and $\mu$ is a measure defined over $(Z_1,\mathcal{F}_1)$, then $f_\sharp \mu$ is a measure over $(Z_2,\mathcal{F}_2)$ defined as 
\[f_\sharp \mu(E) := \mu(f^{-1}(E)) \quad \forall E\in \mathcal{F}_2.\]

Notice that, together with what we introduced §\ref{subsub: space of curve}, we fixed canonical topologies over $\PP(C([a,b],X))$ and $C([a,b],\PP(X))$.

The following measurability result will be useful in the following, and for its proof we refer for example to \cite[Appendix D]{pinzisavare2025}.

\begin{lemma}\label{lemma: meas results}
    Let $g:X \to [0,+\infty]$ and $f:X \to \R^n$ be Borel measurable maps. Then:
    \begin{itemize}
        \item[(i)] the map \[G: \mathcal{M}_+(X)\to[0,+\infty], \quad G(\mu):= \int_X g \ d\mu\]
     is Borel. In particular, the set $\{\mu\in \mathcal{M}_+(X) \ : \ \int_X g d\mu <+\infty\}$ is Borel measurable;
        \item[(ii)] for any $p\geq 1$, the set 
        \[\{(\mu,\nu) \in \mathcal{M}_+(X)\times \mathcal{M}(X;\R^n) \ : \ f\in L^p(\mu), \ \nu = f \mu\}\]
    is Borel, considering the product topology over $\mathcal{M}_+(X)\times \mathcal{M}(X;\R^n)$.
    \end{itemize}
\end{lemma}

\subsubsection{\texorpdfstring{$1$}{}-Wasserstein distance} Let $d$ be a bounded distance that induces the prescribed topology over $X$. The associated $1$-Wasserstein distance, denoted $W_1$, over the space of probability measures is defined as 
    \[W_{1,d} (\mu_1,\mu_2):= \inf \left\{ 
    \int_{X\times X} d(x_1,x_2) d\pi(x_1,x_2) \ : \ \pi \in \Pi(\mu_1,\mu_2) \right\} \quad \text{ for all } \mu_1,\mu_2 \in \PP(X), \]
where $\Pi(\mu_1,\mu_2)$ is the collection of all the transport plans $\pi$ between $\mu_1$ and $\mu_2$, i.e. all the probability measures $\pi \in \PP(X\times X)$ satisfying $\pi(A\times X) = \mu_1(A)$ and $\pi(X\times B) = \mu_2(B)$ for all $A,B \in \mathcal{B}(X)$. If the distance considered is not bounded, we can always replace it by $\hat{d} = d\wedge 1$, that is still a distance that induces the same topology.

The $1$-Wasserstein distance associated with the truncated distance $\hat{d}$, i.e. $W_{1,\hat{d}}$, induces the narrow topology over $\PP(X)$. 

Another important feature of the $1$-Wasserstein distance is its dual formulation: for all $\mu_1,\mu_2 \in \PP(\R^d)$, it holds
\begin{equation}
    W_{1,\hat{d}}(\mu_1,\mu_2)= \sup\left\{ \int_X \phi d\mu_1 - \int_X \phi d \mu_2 \ : \ \phi \in C_b(X) , \ \phi \text{ 1-Lipschitz w.r.t. }\hat{d} \right\}.
\end{equation}
We will make use of this formulation to define different distances inducing the narrow topology over $\PP(\R^d)$ in §\ref{sec: metric}.

\subsubsection{Space of random measures}\label{subsub: RM}
An important space in this paper is the one of random measures over $\R^d$, i.e. the space of probability measures over probability measures $\PP(\PP(\R^d))$. Since $\PP(\R^d)$ with the narrow topology is Polish, we may endow $\PP(\PP(\R^d))$ with the topology described in the previous subsection, that we will call narrow over narrow topology. Moreover, given any bounded distance $D$ over $\PP(\R^d)$ that induces the narrow topology, we know that $W_{1,D}$ induces the narrow over narrow topology.

\subsection{Martingale problem, KFP equation and superposition principle in \texorpdfstring{$\R^d$}{}}\label{sec: fin dim equations}
In this subsection, working in the Euclidean space, we recall the definition of the martingale problem associated with a second order operator $L$, introduced in \cite{stroock1997multidimensional}, the associated Kolmogorov-Fokker-Planck equation on probability measures, and the superposition principle that links them, for which we rely on \cite{figalli2008existence}, \cite{trevisan2016well}.

Here, we will deal with Borel functions
\begin{equation}
    a:[0,T]\times \R^d \to \operatorname{Sym}_+(\R^{d\times d}), \quad b:[0,T]\times \R^d \to \R^d,
\end{equation}
where $\operatorname{Sym}_+(\R^{d\times d})$ are the positive semidefinite square matrices of size $d$. To them, we can associate for any $t\in[0,T]$ the differential operator
\begin{equation}\label{eq: fin dim operator}
L_{b_t,a
_t}\phi(x) := b_t(x)\cdot \nabla\phi(x) + \frac{1}{2}a_t(x):\nabla^2\phi(x) \quad \forall \phi \in C_b^2(\R^d),
\end{equation}
where $A:B:=\sum_{i,j}A_{ij}B_{ij}$ is the scalar product between matrices. Then, we can define the KFP equations and the martingale problem associated to the operator $L$.

\begin{df}\label{def: KFP fin dim}
    Let $\boldsymbol{\mu} = (\mu_t)_{t\in [0,T]}\in C_T(\PP(\R^d))$.
    We say that the Kolmogorov-Fokker-Planck equation $\partial_t\mu_t = L_{b_t,a_t}^*\mu_t$ is satisfied if
    \begin{equation}\label{eq: int kfp}
        \int_0^T \int_{B_R} |b_t(x)| + |a_t(x)| d\mu_t(x) dt<+\infty  \quad \forall R>0,
    \end{equation}
    and for all $\xi\in C_c^1((0,T))$ and $\phi \in C_c^2(\R^d)$ it holds
    \begin{equation}
        \int_0^T \xi'(t) \int_{\R^d} \phi(x)d\mu_t(x)dt = -\int_0^T \xi(t)\int_{\R^d} L_{b_t,a_t}\phi(x) d\mu_t(x) dt.
    \end{equation}
\end{df}

Notice that, as shown in \cite[Remark 2.3]{trevisan2016well}, it is not restrictive to assume that the curve $\boldsymbol{\mu}$ is narrowly continuous, and it actually satisfies
\begin{equation}\label{eq: 2.9}
    \int_s^t \hspace{-0.15cm} \int \partial_t f(r,x) d\mu_r(x) dr = \int f(t,x) d\mu_t(x) - \int f(s,x) d\mu_s(x) - \int_s^t \hspace{-0.15cm} \int L_{b_r,a_r} f(r,x) d\mu_r(x) dr,
\end{equation}
for all $0\leq s < t \leq T$ and $f \in C_c^{1,2}([0,T]\times \R^d)$. 

\begin{df}\label{def: MP solution fin dim}
    Let $\lambda \in \PP(C_T(\R^d))$. We say that it is a solution to the martingale problem associated to $L$, if 
    \begin{equation}\label{eq: integr MP def}
        \int_{C_T(\R^d)} \int_0^T \big(|b_t(\gamma_t)| + |a_t(\gamma_t)|\big)\mathds{1}_{B_R}(\gamma_t)\, dt d\lambda(\gamma) <+\infty \quad \forall R>0,
    \end{equation}
    and for all $\xi\in C_c^1(0,T)$ and $\phi \in C_c^2(\R^d)$ it holds
    \begin{equation}\label{eq: MP cond def}
        [0,T]\ni t \mapsto X^{\xi\phi}_t(\gamma) := \xi(t)\phi(\gamma_t) -\int_0^t \xi'(r)\phi(\gamma_r) + \xi(r)L_r\phi(\gamma_r) dr
    \end{equation}
    is a martingale in the filtered space $\big(C_T(\R^d) , (\mathcal{F}_t)_{t\in[0,T]}, \lambda\big)$, where $\mathcal{F}_t = \sigma (e_r : r\in[0,t])$.
\end{df}

Again, using density of the span of separated variables functions, it is not difficult to see that if $\lambda$ satisfies \eqref{eq: integr MP def} and \eqref{eq: MP cond def}, then it holds that 
\begin{equation}\label{eq: martingale cond with f(t,x)}
[0,T]\ni t \mapsto X_t^{f}(\gamma):= f(t,\gamma_t) - f(0,\gamma_0)-  \int_0^t \partial_t f(r,\gamma_r) + L_{b_r,a_r} f(r,\gamma_r) dr
\end{equation}
is a martingale for each $f\in C_c^{1,2}([0,T]\times \R^d)$. 

\begin{lemma}\label{lemma: equiv of sigma-algebras}
    Let $(X,\tau)$ be a Polish space and consider
    \[\e_s:C_T(X) \to X, \ \  \e_s(\gamma) := \gamma(s), \quad |_{[0,s]}:C_T(X) \to C_s(X), \ \ |_{[0,s]}(\gamma) = \gamma|_{[0,s]}.\]
    Endow $C_t(X)$ with its natural topology and the associated Borel $\sigma$-algebra for any $t\in [0,T]$. Then, for all $t\in [0,T]$, the smallest $\sigma$-algebra on $C_T(X)$ that makes measurable $\e_s$ for all $s\in [0,t]$ coincides with the smallest $\sigma$-algebra that makes measurable $|_{[0,t]}$. Such $\sigma$-algebra will be indicated with $\mathcal{F}_t$, and the collection $(\mathcal{F}_t)_{t\in [0,T]}$ is commonly called the natural filtration of $C_T(X)$.
\end{lemma}

\begin{proof}
    Let $\mathcal{F}_t:=\sigma(\e_s  :  s\in [0,t])$ and $\tilde{\mathcal{F}}_{t} = \sigma(|_{[0,t]})$. The inclusion $\mathcal{F}_t\subseteq \tilde{\mathcal{F}}_{t}$ is immediate, since $\e_s = \e_s^{(t)}\circ |_{[0,t]}$ for all $s\in [0,t]$, where $\e_s^{(t)}:C_t(X) \to X$ is the evaluation at time $s\leq t$ for curves defined up to time $t$.
    \\
    Regarding the other inclusion, endow $X$ with a distance $d$ that makes it complete and separable, and $C_t(X)$ with the sup distance, and notice that for all balls $B_r(\Tilde{\gamma})$, with $\Tilde{\gamma}\in C_t(X)$, it holds
    \begin{align*}
        |_{[0,t]}^{-1}& (B_r(\Tilde{\gamma})) = \{ \gamma \in C_T(X) \ : \ \max_{s\in [0,t]}d(\gamma(s), \Tilde{\gamma}(s)) < r \}
        \\
        =&\bigcup_{k\geq 1} \{ \gamma \in C_T(X) \ : \ \max_{s\in [0,t]} d(\gamma(s), \Tilde{\gamma}(s)) \leq r- 1/k \}
        \\
        = & \bigcup_{k\geq 1} \{ \gamma \in C_T(X) \ : \ \sup_{q\in [0,t]\cap\mathbb{Q}}d(\gamma(q), \Tilde{\gamma}(q)) \leq r- 1/k \}
        \\
        = & \bigcup_{k\geq 1} \bigcap_{q\in [0,t]\cap \mathbb{Q}}\{ \gamma \in C_T(X) \ : \ d(\gamma(q), \Tilde{\gamma}(q)) \leq r- 1/k \}
        \\
        = & \bigcup_{k\geq 1} \bigcap_{q\in [0,t]\cap \mathbb{Q}}\{ \gamma \in C_T(X) \ : \ \e_q(\gamma) \in \overline{B}_{r-1/k}(\Tilde{\gamma}(q))\} \in \mathcal{F}_t.
    \end{align*}
    Since $C_t(X)$ is separable (see \cite[Theorem 2.4.3]{srivastava2008course}), it has a countable basis for the topology made of balls, and then we can conclude that $|_{[0,t]}^{-1}(A) \in \mathcal{F}_t$ for all $A\subset C_t(X)$ open. Then $|_{[0,t]}$ is $\mathcal{F}_t$-measurable, i.e. $\tilde{\mathcal{F}}_{t}\subseteq \mathcal{F}_t$.
\end{proof}


In this work, we will always deal with solutions of some martingale problem, but it is important to recall their connection with solutions to stochastic differential equations. Let $\sigma:[0,T]\times \R^d \to \R^{d\times m}$, $m \geq 1$, be such that $a = \sigma\sigma^\top$. Then $\lambda$ is a martingale solution associated with $L_{b_t,a_t}$ if and only if there exists a probability space $(\Omega,\mathcal{F},\mathbb{P})$ and two continuous processes $(X_t)_{t\in[0,T]}$, $(B_t)_{t\in[0,T]}$ defined on it, and respectively taking values in $
\R^d$ and $\R^m$, such that $B$ is a Brownian motion, the law of $X$ is $\lambda$ and it satisfies 
\begin{equation}\label{eq: SDE}
    dX_t = b_t(X_t) dt + \sigma_t(X_t) dB_t \ \  \text{and} \ \ \  \forall R>0 \ \ \ \mathbb{E}\left[\int_0^T \hspace{-0.25cm}\big( |b_t(X_t)| + |a_t(X_t)| \big)\mathds{1}_{B_R}(X_t) dt\right]<+\infty .
\end{equation}
In particular, uniqueness of solutions of the martingale problem is equivalent to the uniqueness in law of solutions to \eqref{eq: SDE}.

In \cite{figalli2008existence, trevisan2016well, bogachev2021ambrosio}, they proved the so-called \textit{superposition principle}, linking solutions to the martingale problem with solutions of the Kolmogorov-Fokker-Planck equation. Here we state its more general version, specifically referring to \cite[Theorem 1.1]{bogachev2021ambrosio}.

\begin{teorema}\label{thm: fin dim superposition}
    Let $b:[0,T]\times \R^d \to \R^d$ and $a:[0,T]\times\R^d \to \operatorname{Sym}_+(\R^{d\times d})$ be Borel functions, and $L_t:= L_{b_t,a_t}$ their associated operator as in \eqref{eq: fin dim operator}. Let $(\mu_t)_{t\in [0,T]} \in C_T(\PP(\R^d))$ be a solution to $\partial_t\mu_t = L_t^* \mu_t$ satisfying, together with \eqref{eq: int kfp}, the additional integrability assumption
    \begin{equation}\label{eq: gen int ass superposition}
        \int_0^T \int_{\R^d} \frac{|a_t(x)| + |\langle b_t(x),x\rangle|}{1+|x|^2} d\mu_t(x)dt<+\infty.
    \end{equation} 
    Then there exists $\lambda \in \PP(C_T(\R^d))$ solution of the martingale problem associated with $L$, such that $\mu_t = (\e_t)_\sharp \lambda$.
    \\
    Conversely, if $\lambda \in \PP(C_T(\R^d))$ is a solution of the martingale problem, then $\mu^\lambda_t:= (\e_t)_\sharp \lambda$ is a solution for the Kolmogorov-Fokker-Planck equation.
\end{teorema}

\noindent Notice that \eqref{eq: gen int ass superposition} is satisfied if the following stronger assumption holds:
\begin{equation}\label{eq: int cond superposition}
    \int_0^T \int_{\R^d} \frac{|b_t(x)|}{1+|x|} + \frac{|a_t(x)|}{1+|x|^2} d\mu_t(x)dt <+\infty.
\end{equation}

\section{Integral metrics}\label{sec: metric}
In this section, we introduce new integral metrics (see \cite{zolotarev1984probability}) over the space of probability measures. They seem to be quite natural to work in our setting, but we think they could have independent interest. Anyway, we didn't find any previous literature about them.

We will use the notation $C_0^\ell(\R^d)$ to denote functions $\phi \in C^\ell(\R^d)$ such that $\phi$ and all its derivatives of order at most $\ell$ tend to $0$ as $|x|\to +\infty$. Moreover, \[\|\phi\|_{C^\ell}:= \max_{j=0,\dots,\ell} \|\nabla^{\otimes_j}\phi\|_{\infty} = \max_{j=0,\dots,\ell} \sup_{x\in \R^d} |\nabla^{\otimes_j}\phi(x)|,\]
where $|\nabla^{\otimes_j}\phi(x)|$ is the $2$-norm of all the entries of the tensor. 

\begin{df}[Integral metric over positive measures]
    Let $\ell\geq 1$ be an integer number. For all $\mu,\nu\in \mathcal{M}_+(\R^d)$ we define
    \begin{equation}
        D_{C_0^\ell}(\mu,\nu):= \sup\left\{ \left|\int_{\R^d} \phi d\mu - \int_{\R^d} \phi d\nu\right| \ : \ \phi\in C_0^\ell(\R^d), \ \|\phi\|_{C^\ell}\leq 1 \right\}.
    \end{equation}
    Shortly, we will more often use the notation $D_\ell$ for $D_{C_0^\ell}$.
\end{df}

\begin{prop}\label{prop: smooth metrics}
    The following properties hold for any $\ell\geq 1$:
    \begin{enumerate}
        \item $D_\ell$ is a distance over $\mathcal{M}_+(\R^d)$;
        \item $D_\ell$ is a complete metric that metrizes narrow convergence. In particular, $(\PP(\R^d),D_\ell)$ is a complete and separable metric space;
        \item there exists a countable subset $\mathcal{S}_\ell \subset C_c^\ell(\R^d)$ satisfying $\|\phi\|_{C^\ell}\leq 1$ for any $\phi \in \mathcal{S}_\ell$, which is dense in the unit ball of $(C_0^\ell, \|\cdot\|_{C^\ell})$ and 
        \begin{equation}
            D_{C_0^\ell}(\mu,\nu) = \sup\left\{ \left|\int_{\R^d} \phi d\mu - \int_{\R^d} \phi d\nu\right| \ : \ \phi\in \mathcal{S}_\ell \right\}.
        \end{equation}
    \end{enumerate}
\end{prop}

\begin{proof}
    (1) The symmetry and the triangular inequality are trivial. If $D_\ell(\mu,\nu) = 0$, then $\int\phi d\mu = \int \phi d\nu$ for all $\phi \in C_0^\ell(\R^d)$, which is a determining class and thus gives $\mu=\nu$. 

    (2) Let $\mu_n \in \mathcal{M}_+(\R^d)$ be a Cauchy sequence w.r.t. $D_\ell$. We first show that $\sup_n \mu_n(\R^d) <+\infty$. Indeed, consider a non-negative function $\varphi\in C^\infty([0,+\infty))$ such that $\varphi(r) = 1$ in $[0,1]$, $\varphi(r) = 0$ in $[2,+\infty)$ and that is non-increasing. Then, consider $\phi_R(x):= \varphi(|x|/R)$. By monotone convergence theorem, for every $\mu\in \mathcal{M}_+(\R^d)$ it holds
    \[\lim_{R\to +\infty} \int_{\R^d} \phi_R(x) d\mu = \mu(\R^d),\]
    and moreover, $\|\phi_R\|_{C^\ell}\leq 1$ for $R>1$ sufficiently large. Then, considering $\Bar{n}\in \N$ such that for all $n\geq \bar{n}$, $D_\ell(\mu_{n},\mu_{\Bar{n}}) \leq 1$, we have 
    \[ \mu_n(\R^d) = \lim_{R\to+\infty} \int_{\R^d} \phi_R d\mu_n \leq D_\ell(\mu_n,\mu_{\Bar{n}}) + \int_{\R^d} \phi_R d\mu_{\bar{n}}\leq \mu_{\bar{n}}(\R^d) +1,\]
    for all $n\geq \bar{n}$, which allows to conclude.
    
    Now, we show that the sequence $\mu_n$ is tight, i.e. for all $\varepsilon>0$ there exists $R>0$ such that $\mu_n(B_R^c) \leq \varepsilon$ for any $n\in \N$. Indeed, by contradiction, assume that there exists $\varepsilon>0$ such that for all $j \in \N$ there exists $n(j)\to +\infty$ for which $\mu_{n(j)}(B_j^c)>\varepsilon$. In particular, there exists $R(j)\gg j$ such that $\mu_{n(j)}(B_{R(j)} \setminus B_j) \geq \frac{\varepsilon}{2}$. Now, consider a function $\phi_j \in C_c^\ell(\R^d)$ satisfying 
    \[\phi_j \equiv 1 \text{ on }B_{R(j)} \setminus B_j, \quad \phi_j\equiv 0 \text{ on }B_{j/2}\cup B_{2R(j)}^c, \quad \|\phi_j\|_{C^\ell}\leq 1.\]
    Such functions exist at least for $j\geq m$, where $m$ is big enough and depends only on the required regularity $\ell$, and can be built using a cutoff function as done above (see also the proof of the next proposition). Then, consider $\Bar{n} \in \N$ such that $D_\ell(\mu_n,\mu_{\Bar{n}})\leq \varepsilon/4$ for all $n\geq \Bar{n}$. So,
    \[\frac{\varepsilon}{4}\geq \liminf_{j\to+\infty} D_\ell(\mu_{n(j)},\mu_{\Bar{n}}) \geq \liminf_{j \to +\infty} \int\phi_j d\mu_{n(j)} - \int\phi_j d\mu_{\Bar{n}} \geq \frac{\varepsilon}{2} - \limsup_{j\to+\infty}\int \phi_j d\mu_{\Bar{n}} = \frac{\varepsilon}{2},\]
    obtaining a contradiction. Having tightness, we know that there exists $\mu\in \mathcal{M}_+(\R^d)$ such that $\mu_{n_k} \to \mu$ narrowly, for some subsequence $n_k$. 
    It is easy to see that $\mu_{n_k} \to \mu$ narrowly implies $D_\ell(\mu_{n_k},\mu) \to 0$, giving in particular that the whole sequence $\mu_n$ converges to $\mu$ with respect to $D_\ell$ and narrowly. Indeed, if $\mu_{n_k} \to \mu$ narrowly, we know that $\mu_{n_k} \to \mu$ in the \textit{bounded Lipschitz distance} $D_{\operatorname{BL}}$, see e.g. \cite[Section 8.3]{bogachev2007measure}. By its definition we then have
    \begin{align*}
        D_1&(\mu_{n_k},\mu)  = \sup_{\|\phi\|_{C_0^1}\leq 1} \left|\int \phi \, d \mu_{n_k} - \int \phi \, d\mu\right| \\ &\leq \sup \left\{ \int\phi \, d(\mu_{n_k}-\mu)  : \phi\in \operatorname{Lip}_b, \ \|\phi\|_\infty \leq 1, \ \operatorname{LIP}(\phi)\leq 1 \right\}
         = D_{\operatorname{BL}}(\mu_{n_k},\mu) \to 0, 
    \end{align*}
    where $\operatorname{Lip}_b$ is the set of bounded and Lipschitz functions of $\R^d$ and $\operatorname{LIP}(\phi)$ is the global Lipschitz constant of $\phi$.
    Then, notice that for any $\ell>1$ it holds
    \[D_\ell(\mu_{n_k},\mu) \leq D_1(\mu_{n_k},\mu).\]
    In particular, we proved that $\mu_n \to \mu$ narrowly if and only if $D_\ell(\mu_n,\mu) \to 0$.

    (3) It suffices to notice that the Banach space $\big(C_0^\ell(\R^d),\|\cdot\|_{C^\ell}\big)$ is separable and $C_c^\ell(\R^d)$ is dense in it.
\end{proof}

We can also consider an integral metric adding stronger conditions on the test functions. We define it just for the case $\ell = 2$, but it can be easily generalized. Let $C^2_{0,w}(\R^d)$ be the space of \textit{weighted $C^2$-functions}, defined as the closure of $C_c^2(\R^d)$ with respect to the norm
\begin{equation}
    \|\phi\|_{C_{0,w}^2}:= \|\phi\|_\infty + \|(1+|\cdot|)\nabla \phi \|_\infty + \|(1+|\cdot|^2)\nablasquare\phi\|_\infty.
\end{equation}
The space $C_{0,w}^2(\R^d)$ is Banach and separable. We can then consider the \textit{weighted integral metric} defined by 
\begin{equation}
\begin{aligned}
    D_{2,w}(\mu,\nu) = & D_{C_{0,w}^2}(\mu,\nu):= \sup\left\{ \left|\int_{\R^d} \phi d\mu - \int_{\R^d} \phi d\nu\right| \ : \ \phi\in C_{0,w}^2(\R^d), \ \|\phi\|_{C^2_{0,w}}\leq 1 \right\}
    \\
    = & \sup\left\{ \left|\int_{\R^d} \phi d\mu - \int_{\R^d} \phi d\nu\right| \ : \ \phi\in C_{c}^2(\R^d), \ \|\phi\|_{C^2_{0,w}}\leq 1 \right\}.
\end{aligned}
\end{equation}

The analogue of Proposition \ref{prop: smooth metrics} can be proved for this distance.

\begin{prop}\label{prop: smooth weighted metric}
    The following hold:
    \begin{enumerate}
        \item $D_{2,w}$ is a distance over $\mathcal{M}_+(\R^d)$;
        \item $D_{2,w}$ is a complete metric that metrizes narrow convergence. In particular, $(\PP(\R^d),D_{2,w})$ is a complete and separable metric space;
        \item there exists a countable subset $\mathcal{S}_{2,w} \subset C_c^2(\R^d)$ dense in the unit ball of $(C_{0,w}^{2}, \| \cdot \|_{C_{0,w}^2})$, such that 
        \begin{equation}
            D_{2,w}(\mu,\nu) = \sup\left\{ \left|\int_{\R^d} \phi d\mu - \int_{\R^d} \phi d\nu\right| \ : \ \phi\in \mathcal{S}_{2,w} \right\}.
        \end{equation}
    \end{enumerate}
\end{prop}

\begin{proof}
    The proof of Claims (1) and (3) is the same as the one of Proposition \ref{prop: smooth metrics}. The proof of Claim (2) follows its same line, but we must be careful in the choice of the test function considered to obtain a contradiction. So, as before consider $\mu_n\in \mathcal{M}_+(\R^d)$ a Cauchy sequence with respect to $D_{2,w}$.

    The proof that the masses are uniformly bounded works using the same cutoff functions built in the previous proof. Indeed, computing its gradient and its Hessian, it is easy to see that $\|\phi_R\|_{C^2_{0,w}}$ is bounded independently of $R>1$, and thus, up to rescale it, we can use them as test functions for $D_{2,w}$, concluding as before.
    
    Arguing as in the previous proof, we are done if we prove that the sequence is tight. By contradiction, assume that there exists $\varepsilon>0$ such that for all $j\in \N$ there exists $n(j) \to +\infty$ for which $\mu_{n(j)}(B_j^c)>\varepsilon$. In particular, there exists $R(j)\gg j$ such that $\mu_{n(j)}\big( B_{R(j)} \setminus B_j \big) \geq \frac{\varepsilon}{2}$ for all $j\in \N$. 
    \\
    Now, consider $\rho\in C^2([0,2])$ such that $\rho(r)=0$ for all $r\in[0,1]$, $\rho(r)= 1$ in a neighborhood of $r=2$, $\rho$ non-decreasing. We use it to build the following test functions $\phi_j \in C_c^2(\R^d)$:
    \begin{equation}
        \phi_j(x)= \rho_j(|x|):= 
        \begin{cases}
            \rho\left(\frac{2|x|}{j}\right) \quad & \text{if }|x|\leq j,
            \\
            1 &\text{if }|x|\in (j,R(j))
            \\
            \rho\left(3- \frac{|x|}{R(j)}\right) \quad & \text{if }|x|\in[R(j),3R(j))  
            \\
            0 & \text{if }|x|\geq 3R(j).
        \end{cases}
    \end{equation} 
    Computing $\rho_j'(r)$ and $\rho_j''(r)$, one can show that 
    \[|\rho_j'(r)|\leq \frac{3\|\rho'\|_\infty}{r}\mathds{1}_{(\frac{j}{2},+\infty)}(r), \quad |\rho_j''(r)|\leq \frac{9\|\rho''\|_\infty}{r^2}\mathds{1}_{(\frac{j}{2},+\infty)}(r). \]
    This takes us to the estimates on $\phi_j$:
    \[|\nabla \phi_j(x)| \leq \frac{3\|\rho'\|_\infty}{|x|}\mathds{1}_{(\frac{j}{2},+\infty)}(|x|), \quad |\nablasquare\phi_j(x)|\leq \frac{3\sqrt{d-1}\|\rho'\|_\infty + 9\|\rho''\|_\infty}{|x|^2} \mathds{1}_{(\frac{j}{2},+\infty)}(|x|).\]
    
    Fixing then $C:= 1+ 18\|\rho''\|_\infty + 6(1+\sqrt{d-1})\|\rho'\|_\infty$, independent of $j$, and considering the test functions $f_j(x):= \frac{\phi_j(x)}{C}$, it holds that $f_j$ is a competitor in the supremum for $D_{2,w}$, satisfying $f_j(x) = \frac{1}{C}$ for all $|x|\in (j,R(j))$, $0\leq f_j \leq \frac{1}{C}$ and $f_j(x) = 0$ for all $|x|\leq \frac{j}{2}$. Fixing then $\Bar{n}\in \N$ such that $D_{2,w}(\mu_n,\mu_{\Bar{n}})\leq \frac{\varepsilon}{4C}$ for all $n\geq \Bar{n}$, we reach a contradiction:
    \[\frac{\varepsilon}{4C}\geq \liminf_{j\to+\infty} D_{2,w}(\mu_{n(j)},\mu_{\Bar{n}}) \geq \liminf_{j \to +\infty} \int f_j d\mu_{n(j)} - \int f_j d\mu_{\Bar{n}} \geq \frac{\varepsilon}{2C} - \limsup_{j\to+\infty}\int f_j d\mu_{\Bar{n}} = \frac{\varepsilon}{2C}.\]
\end{proof}

Thanks to these distances, we introduce a natural correspondence between the probability measures and the space $\R^\infty$, for which we recall that there are two natural metrics to consider: 
\begin{itemize}
    \item the metric \begin{equation}\label{eq: D_infty distance}
    D_\infty(x,y):= \sup_{n\in \N} |x_n - y_n| \wedge 2,
    \end{equation}
    inducing the uniform convergence (the choice of the constant $2$ will be clear later);
    \item the topology $\tau_w$ induced by the element-wise convergence, i.e. $x\to y$ if $x_n \to y_n $ for all $n\in \N$, induced by the distance 
    \begin{equation}\label{eq: element-wise distance}
        d_\infty(x,y) := \sum_{n\in \N} \frac{|x_n-y_n|\wedge 1}{2^n}.
    \end{equation}
    The topological space $(\R^\infty,\tau_w)$ is Polish.
\end{itemize}
Such a space is strictly related to the space of probability measures when endowed with the smooth metrics $D_\ell$: fix $\ell\geq 1$ and say that $\mathcal{S}_\ell = \{\varphi_{1,\ell},\varphi_{2,\ell},\dots\}$ as the countable subset given by Proposition \eqref{prop: smooth metrics}, (3). Then we can define
\begin{equation}\label{eq: def iota}
\begin{aligned}
    \iota_\ell: \PP(\R^d) & \to \R^{\infty}
    \\
    \mu & \mapsto (L_{\varphi_{1,\ell}}(\mu),L_{\varphi_{2,\ell}}(\mu), \dots)
\end{aligned}
\end{equation}
We are mainly interested in the case $\ell = 2$ and this construction will be fundamental for the proof of our main theorem, and for this reason, we introduce also the map 
\begin{equation}\label{eq: weighted map R^infty}
    \begin{aligned}
    \iota_{2,w}: \PP(\R^d) & \to \R^{\infty}
    \\
    \mu & \mapsto (L_{\varphi_{1}^{(w)}}(\mu),L_{\varphi_{2}^{(w)}}(\mu), \dots),
\end{aligned}
\end{equation}
where $\mathcal{S}_{2,w}= \{\varphi_{1}^{(w)}, \varphi_{2}^{(w)},\dots\}$ is the countable set of test functions introduced in Proposition \ref{prop: smooth weighted metric}. 

An easy but useful corollary is the following.

\begin{co}\label{cor: closed and isometry}
    The map $\iota_\ell$ is continuous and injective between $(\PP(\R^d),D_\ell)$ and $(\R^\infty, d_\infty)$. Moreover, the set $\iota_\ell(\PP(\R^d))$ is closed in $(\R^\infty,D_\infty)$, and $\iota_\ell$ is an isometry between $(\PP(\R^d),D_\ell)$ and $(\iota_\ell(\PP(\R^d)), D_\infty)$. The same holds for $\iota_{2,w}$ and $D_{2,w}$.
\end{co}

\begin{proof}
    We first prove it is an isometry. Consider $\mu,\nu \in \PP(\R^d)$. Then, denoting $\mathcal{S}_\ell= \{\varphi_1,\varphi_2, \dots,\}$, thanks to Proposition \ref{prop: smooth metrics},
    \[D_\infty(\iota_\ell(\mu), \iota_\ell(\nu)) = \sup_{n\in \N} \left| \int_{\R^d} \varphi_n d\mu - \int_{\R^d} \varphi_n d\nu \right| \wedge 2 = D_\ell(\mu,\nu)\wedge 2 = D_\ell(\mu,\nu),\]
    where the last equality follows from the fact that $\|\varphi_n\|\leq 1$ for every $n\in \N$.
    Since $(\PP(\R^d),D_\ell)$ is complete, its isometric image is complete, which concludes the proof since complete subspaces of a metric space are closed.
    The proof is the same for the map $\iota_{2,w}$, relying on Proposition \ref{prop: smooth weighted metric}.
\end{proof}

\section{Equation on random measures and nested superposition principle}\label{sec: nd equations}
In this section, we introduce an evolution equation for random measures associated with the Borel functions
\begin{equation}\label{eq: non-local vf}
    a:[0,T]\times \R^d \times \PP(\R^d) \to \operatorname{Sym}_+(\R^{d\times d}), \quad b:[0,T]\times \R^d \times \PP(\R^d) \to \R^d.
\end{equation}

According to the Kolmogorov-Fokker-Planck equation and the martingale problem presented in Section \ref{sec: fin dim equations}, we introduce the operator
\begin{equation}\label{eq: non-local operator 1}
    L_{b,a} \phi(t,x,\mu) := \frac{1}{2}a_t(x,\mu):\nabla^{\otimes_2}\phi(x) + b_t(x,\mu)\cdot \nabla \phi(x),
\end{equation}
for any $\phi \in C_b^2(\R^d)$. It is also useful to define the following objects: given $\boldsymbol{\mu}  = (\mu_t)_{t\in[0,T]} \in C_T(\PP(\R^d))$ and $\lambda \in \PP(C_T(\R^d))$, for any $\phi \in C_b^2(\R^d)$ let 
\begin{equation}\label{eq: non-local operators 2}
\begin{aligned}
    & L^{\boldsymbol{\mu}}_{b_t,a_t} \phi(x) := \frac{1}{2}a_t(x,\mu_t):\nabla^{\otimes_2}\phi(x) + b_t(x,\mu_t)\cdot \nabla \phi(x),
    \\
    & L^\lambda_{b_t,a_t}\phi(x) := \frac{1}{2}a_t(x,(\e_t)_\sharp \lambda):\nabla^{\otimes_2}\phi(x) + b_t(x,(\e_t)_\sharp \lambda)\cdot \nabla \phi(x).
\end{aligned}
\end{equation}

    \noindent Notice that, for any $\phi \in C^2_b$, it holds
    \begin{equation}\label{eq: estimate for the operator}
    |L_{b,a}\phi(t,x,\mu)| \leq \frac{1}{2}|\nablasquare\phi||a(t,x,\mu)| + |\nabla\phi| |b(t,x,\mu)|,\end{equation}
    and similarly for the operator written in the forms $L_{b_t,a_t}^{\boldsymbol{\mu}}\phi$ and $L_{b_t,a_t}^\lambda \phi$.

    The non-local nature of the operator $L$ leads to an operator that acts on cylinder functions.

    \begin{df}[Cylinder functions]\label{def: cyl functions}
    Let $h,\ell \geq 0$. A functional $F: \PP(\R^d) \to \R$ is said to belong to $\operatorname{Cyl}_c^{h,\ell}(\PP(\R^d))$ if there exists $k\in \N$, $\Psi \in C_c^h(\R^k)$ and $\Phi = (\phi_1,\dots,\phi_k) \in C_c^\ell(\R^d,\R^k)$ such that 
    \begin{equation}\label{eq: cyl functions}
        F(\mu) = \Psi\left( L_\Phi(\mu) \right), \quad L_\Phi(\mu) =  \big(L_{\phi_1}(\mu) ,\dots , L_{\phi_k}(\mu)\big), \quad  L_{\phi_i}(\mu) := \int_{\R^d} \phi_i(x) d\mu(x).
    \end{equation}
    If $\Phi \in C_b^\ell(\R^d,\R^k)$ and $\Psi \in C_b^h(\R^k)$, then we say that $F\in \operatorname{Cyl}_b^{h,\ell}(\PP(\R^d))$.
\end{df}
Notice that, if $\Psi \circ L_\Phi \in \operatorname{Cyl}_c^{h,\ell}(\PP(\R^d))$, we can also consider $\Psi \in C_b^h(\R^k)$. Moreover, $\operatorname{Cyl}_c^{h,\ell}(\PP(\R^d))\subset\operatorname{Cyl}_b^{h,\ell}(\PP(\R^d))$.
    
    Now, for all $t\in[0,
    T]$ we define the operator $\mathcal{K}_{b_t,a_t}$ naturally acting on cylinder functions $\operatorname{Cyl}_b^{1,2}(\PP(\R^d))$ as
    \begin{equation}
    \begin{aligned}
        \mathcal{K}_{b_t,a_t} F(x,\mu) := & \sum_{i=1}^k \partial_i \Psi(L_\Phi(\mu)) L_{b,a}\phi_i(t,x,\mu)
        \\
        = & \sum_{i=1}^k \partial_i \Psi(L_\Phi(\mu))\Big( b_t(x,\mu)\cdot \nabla\phi_i(x) + \frac{1}{2}a_t(x,\mu):\nablasquare \phi_i(x) \Big),
    \end{aligned}
    \end{equation}
    for all $F = \Psi\circ L_{\boldsymbol{\Phi}} \in \operatorname{Cyl}_b^{1,2}(\PP(\R^d))$. It is well defined, i.e. it does not depend on the representation chosen for the cylinder function $F$, since it is uniquely determined by the expressions
    \[\sum_{i=1}^k \partial_i \Psi(L_\Phi(\mu))\nabla\phi_i(x) = \nabla_x \left(\frac{d^+}{d\varepsilon}|_{\varepsilon=0} \ F((1-\varepsilon)\mu+ \varepsilon\delta_x)\right),\]
    \[\sum_{i=1}^k \partial_i \Psi(L_\Phi(\mu))\nablasquare\phi_i(x) = \nablasquare_x \left(\frac{d^+}{d\varepsilon}|_{\varepsilon=0} \ F((1-\varepsilon)\mu+ \varepsilon\delta_x)\right).\]

\begin{df}\label{def: sol of nd KFP for rm}
    Let $a$ and $b$ as in \eqref{eq: non-local vf} and $\boldsymbol{M} = (M_t)_{t\in [0,T]} \in C_T(\PP(\PP(\R^d)))$. We say that it solves the equation $\partial_t M_t = \mathcal{K}_{b_t,a_t}^* M_t$ if 
    \begin{equation}\label{eq: int cond for nd KFP}
        \int_0^T \int_{\PP}\int_{B_R} |a(t,x,\mu)| + |b(t,x,\mu)| d\mu(x) dM_t(\mu) dt <+\infty \quad \forall R>0,
    \end{equation}
    and for any $F = \Psi\circ L_{\Phi} \in \operatorname{Cyl}_c^{1,2}(\PP(\R^d))$ and $\xi \in C_c^1(0,T)$ it holds
    \begin{equation}\label{eq: def nd KFP}
        \int_0^T \xi'(t) \int_{\PP} F(\mu) dM_t(\mu) dt = -\int_0^T \xi(t) \int_{\PP}  \int_{\R^d} \mathcal{K}_{b_t,a_t}F(x,\mu) d\mu(x)dM_t(\mu) dt.
    \end{equation}
\end{df}


The assumption that the curve $t\mapsto M_t \in \PP(\PP(\R^d))$ is continuous is not restrictive for our purposes, as we will always assume the stronger integrability condition
\begin{equation}\label{eq: weighted int cond random measures}
    \int_0^T \int_{\PP}\int_{\R^d} \frac{|a(t,x,\mu)|}{1+|x|^2} + \frac{|b(t,x,\mu)|}{1+|x|} d\mu(x) dM_t(\mu) dt <+\infty.
\end{equation}

\begin{lemma}\label{lemma: sol of ndKFP are AC}
    Let $a$ and $b$ be as above and $\boldsymbol{M}=(M_t)_{t\in [0,T]}\subset \PP(\PP(\R^d))$ such that $[0,T]\ni t\mapsto M_t \in \PP(\PP(\R^d))$ is Borel measurable and satisfies 
    \eqref{eq: weighted int cond random measures}
    and \eqref{eq: def nd KFP}. Then, there exists a curve $(\Tilde{M}_t)_{t\in [0,T]} \in AC_T(\PP(\PP(\R^d)), \mathcal{W}_{1,D_{2,w}})$ such that $\Tilde{M}_t = M_t$ for a.e. $t\in [0,T]$. In particular, $(\Tilde{M}_t)_{t\in [0,T]}\in C_T(\PP(\PP(\R^d)))$ and is the unique continuous representative for $\boldsymbol{M}$.
\end{lemma}

\begin{proof}
    For any $F = \Psi \circ L_{\Phi} \in \operatorname{Cyl}_c^{1,2}$, the function $t\mapsto \int F(\mu)dM_t(\mu)$ is in $W^{1,1}(0,T)$, with distributional derivative 
    \[t \mapsto \int_{\PP} \sum_{i=1}^k \partial_i \Psi(L_\Phi(\mu)) \int_{\R^d} L_{b,a}\phi_i(t,x,\mu) d\mu(x)dM_t(\mu) \in L^1(0,T).\]
    Then, there exists $I_F\subset(0,T)$ such that for any $s,t \in I_F$
    \[\int F(\mu)dM_t(\mu) - \int F(\mu)dM_s(\mu) = \int_s^t \int_{\PP} \sum_{i=1}^k \partial_i \Psi(L_\Phi(\mu)) \int_{\R^d} L_{b,a}\phi_i(r,x,\mu) d\mu(x)dM_r(\mu) dr. \]
    Consider now $I:= \bigcap_{F \in \mathcal{C}}I_F$ with $\mathcal{C}$ satisfying the following properties:
    \begin{itemize}
        \item $F = \Psi \circ L_{\Phi} \in \mathcal{C}$ implies that $\Psi \in \mathcal{F}_k$ for some $k\in \N$ countable such that $\Psi:\R^k \to \R$ is continuously differentiable satisfying $\|\Psi\|_{\infty} \leq 1,\llbracket\Psi\rrbracket_{C^1}:= \sup_{y\in \R^k} \sum_{i=1}^k |\partial_i \Psi(y)| \leq 1$ and 
        \[W_{1,D_\infty}(\widetilde{M},\widetilde{N}) = \sup_{k\in \N} \sup_{\Psi \in \mathcal{F}_k} \int\Psi(x_1,\dots,x_k) d\big(\widetilde{M} - \widetilde{N}\big)(\underline{x}).\]
        See \cite[Remark C.7]{pinzisavare2025} for the existence of the $\mathcal{F}_k$'s;
        \item $F = \Psi \circ L_{\Phi}$ implies that $\Phi = (\varphi_1,\dots, \varphi_k)$ for some $k\in \N$ and $\mathcal{S}_{2,w} = \{\varphi_1,\varphi_2,\varphi_3,\dots\}$ introduced in Proposition \ref{prop: smooth weighted metric}, (3). To make the notation clearer, we write $L_{\Phi} = L_{k,\boldsymbol{\varphi}}$.
    \end{itemize}
    These properties imply that $\mathcal{C}$ is countable, so $I\subset(0,T)$ has full Lebesgue measure. In particular, defining $\iota:=\iota_{2,w}$ as in \eqref{eq: weighted map R^infty}, for any $s,t \in I$ it holds
    \[
    \begin{aligned}
    W_{1,D_{2,w}} & (M_s,M_t) = W_{1,D_{\infty}}(\iota_{\sharp }M_s,\iota_\sharp  M_t) = \sup_{k\in \N} \sup_{\Psi \in \mathcal{F}_k}\int\Psi(\pi^k(x_1,\dots,x_k) d\big(\iota_{\sharp }M_s - \iota_{\sharp }M_t\big)(\underline{x})
    \\
    = & 
    \sup_{k\in \N} \sup_{\Psi \in \mathcal{F}_k} \int\Psi(L_{\varphi_1}(\mu),\dots, L_{\varphi_k}(\mu)) d\big(M_s - M_t\big)(\mu)
    \\
    = & \sup_{k\in \N} \sup_{\Psi \in \mathcal{F}_k} \int_s^t \int_{\PP} \sum_{i=1}^k \partial_i \Psi(L_{k,\boldsymbol{\varphi}}(\mu)) \int_{\R^d} L_{b,a}\varphi_i(r,x,\mu) d\mu(x)dM_r(\mu) dr
    \\
    \leq &
    \int_s^t \int_{\PP} \int_{\R^d} \frac{|a(r,x,\mu)|}{1+|x|^2} + \frac{|b(r,x,\mu)|}{1+|x|} d\mu(x)dM_r(\mu) dr,
    \end{aligned}
    \]
    thanks to \eqref{eq: estimate for the operator} and the facts that for all $\Psi \in \mathcal{F}_k$, it holds $\llbracket\Psi \rrbracket_{C^1}\leq 1$ and for all $\varphi_i\in \mathcal{S}_{2,w}$ we have $\|(1+|x|)\nabla \varphi_i\|_\infty\leq 1$ and $\|(1+|x|^2)\nablasquare \varphi_i\|_\infty\leq 1$. In particular, by completeness of $(\PP(\PP(\R^d)), W_{1,D_{2,w}})$, there exists $(\Tilde{M}_t)\in AC_T(\PP(\PP(\R^d)), W_{1,D_{2,w}})$, which satisfies the requirements by construction.
\end{proof}

Let us highlight some features of the operator $\mathcal{K}$: for simplicity we get rid of the time variable, considering two Borel measurable maps $a:\R^d\times \PP(\R^d) \to \operatorname{Sym}_+(\R^{d\times d})$ and $b:\R^d \times \PP(\R^d) \to \R^d$ and the operator 
\[\mathcal{K}_{b,a}F(x,\mu) := \sum_{i=1}^k \partial_i \Psi(L_\Phi(\mu))\Big( b(x,\mu)\cdot \nabla\phi_i(x) + \frac{1}{2}a(x,\mu):\nablasquare \phi_i(x) \Big),\]
for all $F = \Psi\circ L_{\Phi} \in \operatorname{Cyl}_b^{1,2}(\PP(\R^d))$. It is not hard to verify that the operator $\mathcal{K}_{b,a}$ satisfies the Leibniz rule, that is 
\begin{equation}
    \mathcal{K}_{b,a}(FG) = F \mathcal{K}_{b,a} G + G \mathcal{K}_{b,a} F, 
\end{equation}
which is a prerogative of first-order operators. Indeed, we may see the equation $\partial_tM_t = \mathcal{K}_{b_t,a_t}^*M_t$ as a first order equation on the space of probability measures, but to do so it is important to understand that the natural geometry to consider on $\PP(\R^d)$ cannot be the one given by the Wasserstein on Wasserstein distance, but it should be a second order metric that naturally puts probability measures in duality with two times differentiable functions, as it can be seen in the proof of Lemma \ref{lemma: sol of ndKFP are AC} or later in the proof of the superposition result. 

We proceed by showing how this equation on random measures is linked to the KFP equation and the martingale problem introduced in Section \ref{sec: fin dim equations}. In particular, we link the following objects through a superposition principle:
\begin{enumerate}
    \item a curve of random measures $\boldsymbol{M} = (M_t)_{t\in [0,T]} \in C_T(\PP(\PP(\R^d)))$ solution of $\partial_t M_t = \mathcal{K}_{b_t,a_t}^*M_t$, according to Definition \ref{def: sol of nd KFP for rm};
    \item a probability measure $\Lambda \in \PP(C_T(\PP(\R^d)))$ such that 
    \begin{equation}\label{eq: int cond for Lambda nd}
        \int \int_0^T \int \frac{|b(t,x,\mu_t)|}{1+|x|} + \frac{|a(t,x,\mu_t)|}{1+|x|^2} d\mu_t(x) d\Lambda dt (\boldsymbol{\mu}) <+\infty
    \end{equation}
    and $\Lambda$-a.e. $\boldsymbol{\mu} = (\mu_t)_{t\in [0,T]}$ is a solution of $\partial_t \mu_t = (L^{\boldsymbol{\mu}}_{b_t,a_t})^* \mu_t$, according to Definition \ref{def: KFP fin dim};
    \item a random measure over curves $\mathfrak{L}\in \PP(\PP(C_T(\R^d)))$ satisfying 
    \begin{equation}\label{eq: int cond for L nd}
        \int \int \int_0^T \frac{|b(t,\gamma_t,(\e_t)_\sharp \lambda)|}{1+|\gamma_t|} + \frac{|a(t,\gamma_t,(\e_t)_\sharp \lambda)|}{1+|\gamma_t|^2} dt d\lambda(\gamma) d\mathfrak{L}(\lambda)  <+\infty
    \end{equation}
    and $\mathfrak{L}$-a.e. $\lambda\in \PP(C_T(\R^d))$ is a solution to the martingale problem associated to the operator $L^\lambda_{a_t,b_t}$, according to Definition \ref{def: MP solution fin dim}.
\end{enumerate}

\noindent To this aim, it will be crucial to define the following subsets.

\begin{df}\label{def: KFP(b,a) and MP(b,a)}
    Let $b:[0,T]\times \R^d \times \PP(\R^d) \to \R^d$ and $a:[0,T]\times \R^d \times \PP(\R^d) \to \operatorname{Sym}_+(\R^{d\times d})$ be Borel measurable functions. The set $\KFP(b,a) \subset C_T(\PP(\R^d))$ is the subset of the solutions to the non-local Kolmogorov-Fokker-Planck equation, i.e. 
    \begin{equation}\label{eq: KFP(b,a)}
    \begin{aligned}
        \KFP(b,a) := \bigg\{\boldsymbol{\mu} = (\mu_t) \, : \, 
        \int_0^T \hspace{-0.2cm} \int_{\R^d} \frac{|b_t(x,\mu_t)|}{1+|x|} + \frac{|a_t(x,\mu_t)|}{1+|x|^2} d\mu_t(x) dt <+\infty, \ \partial_t \mu_t = (L^{\boldsymbol{\mu}}_{b_t,a_t})^* \mu_t
        \bigg\}.
    \end{aligned}
    \end{equation}
    The set $\MP(b,a) \subset \PP(C_T(\R^d))$ is the subset of solutions to the non-local martingale problem, i.e.
    \begin{equation}\label{eq: MP(b,a)}
    \begin{aligned}
        \MP(b,a) := \bigg\{\lambda \in \PP(C_T(\R^d))  \ : \ 
        & \int\int_0^T \frac{|b_t(\gamma_t,(\e_t)_\sharp \lambda)|}{1+|\gamma_t|} + \frac{|a_t(\gamma_t,(\e_t)_\sharp \lambda)|}{1+|\gamma_t|^2} dtd\lambda(\gamma) <+\infty, \ \\
        & \lambda \text{ is a sol. of the martingale problem associated to }L^\lambda_{b_t,a_t}
        \bigg\}.
    \end{aligned}
    \end{equation}
\end{df}

It is important to stress here the relation between a probability measure $\lambda\in \operatorname{MP}(b,a)$ and solutions of some stochastic differential equation. As pointed out in §\ref{sec: fin dim equations}, if the map $a$ arises from some $\sigma:[0,T]\times \R^d \times \PP(\R^d) \to \R^{d\times m}$, $m\geq 1$, as $a = \sigma\sigma^\top$, then $\lambda \in \operatorname{MP}(b,a)$ if and only if there exists a probability space $(\Omega,\mathcal{F},\mathbb{P})$, an $\R^m$-valued Brownian motion $(B_t)_{t\in[0,T]}$ on it and an $\R^d$-valued continuous process $(X_t)_{t\in[0,T]}$ whose law coincides with $\lambda$, and in particular $(\e_t)_\sharp \lambda = \operatorname{Law}(X_t)$ for all $t\in[0,T]$, so that it satisfies
\begin{equation}\label{eq: sec 4 mckean vlasov}
dX_t = b_t(X_t,\operatorname{Law}(X_t)) dt + \sigma_t(X_t,\operatorname{Law}(X_t))dB_t 
\end{equation}
and 
\begin{equation}\label{eq: int cond mckean vlasov}
    \mathbb{E} \left[ \int_0^T \frac{|b_t(X_t,\operatorname{Law}(X_t))|}{1+|X_t|} + \frac{|a_t(X_t,\operatorname{Law}(X_t))|}{1+|X_t|^2}  \right]<+\infty
\end{equation}

Now, we show that there is a natural hierarchy between the objects presented above. Let us first introduce the following operators:
\begin{equation}\label{eq: maps for marginals}
    \begin{aligned}
        & E:\PP(C_T(\R^d)) \to C_T(\PP(\R^d)), \quad E(\lambda):= ((\e_t)_\sharp \lambda)_{t\in [0,T]},
        \\
        & E_t:\PP(C_T(\R^d)) \to \PP(\R^d), \quad E_t(\lambda) = (\e_t)_\sharp \lambda,
        \\
        & \mathfrak{e}_t: C_T(\PP(\R^d)) \to \PP(\R^d), \quad \mathfrak{e}_t(\boldsymbol{\mu}) = \mu_t.
    \end{aligned}
\end{equation}

\begin{prop}\label{prop: from L to Lambda nd}
    Let $b:[0,T]\times \R^d \times \PP(\R^d) \to \R^d$ and $a:[0,T]\times \R^d \times \PP(\R^d) \to \operatorname{Sym}_+(\R^{d\times d})$ be Borel measurable maps. 
    \\
    The functional $E$ maps $\MP(b,a)$ into $\KFP(b,a)$. In particular, if $\mathfrak{L}\in \PP(\PP(C_T(\R^d)))$ is concentrated over $\MP(b,a)$, then $\Lambda := E_\sharp \mathfrak{L} \in \PP(C_T(\PP(\R^d)))$ is concentrated over $\KFP(b,a)$.
\end{prop}

\begin{proof}
    Let $\lambda \in \MP(b,a)$ and $\mu_t:= (\e_t)_\sharp \lambda$. The integrability condition is clearly satisfied. On the other hand, notice that, directly by definition $L_{b_t,a_t}^\lambda \phi(x) = L_{b_t,a_t}^{\boldsymbol{\mu}}\phi(x)$ for any $x\in \R^d$ and $\phi \in C_c^2(\R^d)$. Then, for any $\xi \in C_c^1(0,T)$ and $\phi \in C_c^2(\R^d)$, for any $t\in[0,T]$ define, as in \eqref{eq: MP cond def},
    \[X_t^{\xi,\phi,\lambda} (\gamma):= \xi(t) \phi(\gamma_t) - \int_0^t \xi'(r)\phi(\gamma_r) + \xi(r)L_{b_r,a_r}^\lambda\phi(\gamma_r) dr.\]
    Since $\lambda\in\MP(b,a)$, $t\mapsto X_t^{\xi,\phi,\lambda}$ is a martingale in the filtered space $(C_T(\R^d), (\mathcal{F}_t)_{t\in [0,T]}, \lambda)$ (see Definition \ref{def: MP solution fin dim}), then
    \begin{align*}
        \int_0^T \int_{\R^d} \xi'(t)\phi(x) + \xi(t)L_{b_t,a_t}^{\boldsymbol{\mu}}\phi(x) d\mu_t(x) dt = \int -X^{\xi,\phi,\lambda}_T(\gamma)d\lambda = \int -X^{\xi,\phi,\lambda}_0(\gamma)d\lambda = 0.
    \end{align*}
\end{proof}

\begin{prop}\label{prop: from Lambda to M nd}
     Let $b:[0,T]\times \R^d \times \PP(\R^d) \to \R^d$ and $a:[0,T]\times \R^d \times \PP(\R^d) \to \operatorname{Sym}_+(\R^{d\times d})$ be Borel measurable functions. 
     \\
     Let $\Lambda \in \PP(C_T(\PP(\R^d)))$ be satisfying \eqref{eq: int cond for Lambda nd} and concentrated over $\KFP(b,a)$. Then the curve of random measures $\boldsymbol{M} \in C_T(\PP(\PP(\R^d)))$, defined by $M_t := (\mathfrak{e}_t)_\sharp \Lambda$, satisfies \eqref{eq: weighted int cond random measures} and solves the equation $\partial_t M_t = \mathcal{K}_{b_t,a_t}^*M_t$, according to Definition \ref{def: sol of nd KFP for rm}.
\end{prop}

\begin{proof}
    First of all, thanks to \eqref{eq: int cond for Lambda nd} and Fubini's theorem, it holds 
    \[
    \begin{aligned}
    \int_0^T \int_{\PP}\int_{\R^d} & \frac{|a(t,x,\mu)|}{1+|x|^2} + \frac{|b(t,x,\mu)|}{1+|x|} d\mu(x) dM_t(\mu) dt =
    \\
    = & \int \int_0^T \int_{\R^d} \frac{|a(t,x,\mu_t)|}{1+|x|^2} + \frac{|b(t,x,\mu_t)|}{1+|x|} d\mu_t(x) dt d\Lambda(\boldsymbol{\mu})<+\infty.
    \end{aligned}
    \]
    Now, let $\xi \in C_c^1(0,T)$ and $F = \Psi \circ L_\Phi \in \operatorname{Cyl}_c^{1,2}(\PP(\R^d))$, then
    \begin{align*}
        \int_0^T & \xi'(t)\int F(\mu)dM_t(\mu)dt = \int_0^T \xi'(t) \int F(\mu_t)d\Lambda(\boldsymbol{\mu}) dt 
        \\
        & = \int \int_0^T \xi'(t)\Psi(L_\Phi(\mu_t))dt d\Lambda(\boldsymbol{\mu})
        \\
        & = - \int \int_0^T \xi(t) \sum_{i=1}^k \partial_i\Psi(L_\Phi(\mu_t)) \frac{d}{dt} \left( \int_{\R^d} \phi_i(x) d\mu_t(x) \right) dt d\Lambda(\boldsymbol{\mu})
        \\
        & = - \int \int_0^T \xi(t)\sum_{i=1}^k \partial_i\Psi(L_\Phi(\mu_t)) \left( \int_{\R^d} L_{b_t,a_t}^{\boldsymbol{\mu}} \phi_i(x) d\mu_t(x)\right) dt d\Lambda(\boldsymbol{\mu})
        \\
        & = - \int \int_0^T \xi(t)\sum_{i=1}^k \partial_i\Psi(L_\Phi(\mu_t)) \left( \int_{\R^d} L_{b,a} \phi_i(t,x,\mu_t) d\mu_t(x)\right) dt d\Lambda(\boldsymbol{\mu})
        \\
        & = - \int_0^T \int \xi(t)\sum_{i=1}^k \partial_i\Psi(L_\Phi(\mu_t)) \left( \int_{\R^d} L_{b,a} \phi_i(t,x,\mu_t) d\mu_t(x)\right) d\Lambda(\boldsymbol{\mu}) dt
        \\
        & = - \int_0^T \int \xi(t)\sum_{i=1}^k \partial_i\Psi(L_\Phi(\mu)) \int_{\R^d} L_{b,a} \phi_i(t,x,\mu) d\mu(x) dM_t(\mu) dt
        \\
        & = - \int_0^T \xi(t) \int \mathcal{K}_{b_t,a_t}F(x,\mu) d\mu(x)dM_t(\mu)dt. 
    \end{align*}
\end{proof}

\begin{co}
Let $b:[0,T]\times \R^d \times \PP(\R^d) \to \R^d$ and $a:[0,T]\times \R^d \times \PP(\R^d) \to \operatorname{Sym}_+(\R^{d\times d})$ be Borel measurable functions. 
     \\
     Let $\mathfrak{L} \in \PP(\PP(C_T(\R^d)))$ be satisfying \eqref{eq: int cond for L nd} and concentrated over $\MP(b,a)$. Then the curve of random measures $\boldsymbol{M} \in C_T(\PP(\PP(\R^d)))$, defined by $M_t := (E_t)_\sharp \mathfrak{L}$, satisfies \eqref{eq: weighted int cond random measures} and solves the equation $\partial_t M_t = \mathcal{K}_{b_t,a_t}^*M_t$, according to Definition \ref{def: sol of nd KFP for rm}.
\end{co}

\begin{proof}
    Noticing that $E_t = \mathfrak{e}_t\circ E$, it follows from Propositions \ref{prop: from L to Lambda nd} and \ref{prop: from Lambda to M nd}.
\end{proof}

To summarize: we proved that there is a natural hierarchy between the objects $\boldsymbol{M}$, $\Lambda$ and $\mathfrak{L}$ listed above, that is $\mathfrak{L} \implies \Lambda \implies \boldsymbol{M}$. In the next subsections, we will show that we can also go in the opposite direction.

\subsection{Superposition principle: from \texorpdfstring{$\boldsymbol{M}$}{} to \texorpdfstring{$\Lambda$}{}}
Here, we show that the result of Proposition \ref{prop: from Lambda to M nd} can be inverted, i.e. given a curve of random measures $\boldsymbol{M} \in C_T(\PP(\PP(\R^d)))$ solving the equation $\partial_t M_t = \mathcal{K}_{b_t,a_t}^*M_t$, we show it is possible to obtain a measure $\Lambda\in \PP(C_T(\PP(\R^d)))$ concentrated over $\KFP(b,a)$ such that $(\mathfrak{e}_t)_\sharp \Lambda = M_t$ for all $t\in[0,T]$. 

\begin{teorema}\label{thm: superposition from M to Lambda}
    Let $\boldsymbol{M}= (M_t)_{t\in [0,T]} \in C_T(\PP(\PP(\R^d)))$. Let $b:[0,T]\times \R^d \times \PP(\R^d) \to \R^d$ and $a:[0,T]\times \R^d \times \PP(\R^d) \to \operatorname{Sym}_+(\R^{d\times d})$ be Borel measurable functions and assume that \eqref{eq: weighted int cond random measures} and \eqref{eq: def nd KFP} are satisfied. 
    Then, there exists (possibly non-unique) $\Lambda \in \PP(C_T(\PP(\R^d)))$ satisfying:
    \begin{enumerate}
        \item $(\mathfrak{e}_t)_\sharp \Lambda = M_t$ for all $t\in [0,T]$;
        \item \[\int \int_0^T \int \frac{|b(t,x,\mu_t)|}{1+|x|} + \frac{|a(t,x,\mu_t)|}{1+|x|^2}d\mu_t(x) dt d\Lambda(\boldsymbol{\mu})<+\infty;\]
        \item $\Lambda$ is concentrated over $\KFP(b,a)$, in particular $\Lambda$-a.e. $\boldsymbol{\mu} = (\mu_t)_{t\in [0,T]}$ solves the Kolmogorov-Fokker-Planck equation $\partial_t \mu_t = (L_{b_t,a_t}^{\boldsymbol{\mu}})^*\mu_t$. 
    \end{enumerate}
\end{teorema}

\begin{proof}
    \textbf{Step 1}: let $\mathcal{S}_{2,w} = \{\varphi_1,\varphi_2,\dots\}$ as in Proposition \ref{prop: smooth weighted metric}, (3) and define $\iota = \iota_{2,w} :\PP(\R^d) \to \R^\infty$ as in \eqref{eq: weighted map R^infty}, keeping in mind that it is an isometry between $(\PP(\R^d),D_{C_{0,w}^2})$ and $(\iota(\PP(\R^d)),D_\infty)$, as shown in Corollary \ref{cor: closed and isometry}. Let $\widetilde{M}_t := \iota_\sharp  M_t$ for any $t\in [0,T]$. Then, $(\widetilde{M}_t)_{t\in [0,T]} \in AC_T(\PP(\R^\infty),W_{1,D_\infty})$, indeed considering any $\Pi_{t,s} \in \Gamma(M_t,M_s)$ optimal for the distance $\mathcal{W}_{1,D_{2,w}}$, then
    \[W_{1,D_\infty}(\widetilde{M}_t,\widetilde{M}_s) \leq \int D_\infty(\underline{x},\underline{y}) d(\iota,\iota)_\sharp \Pi_{t,s}(\underline{x},\underline{y}) = \int D_{2,w}(\mu,\nu) d\Pi_{t,s}(\mu,\nu) = \mathcal{W}_{1,D_{2,w}}(M_t,M_s), \]
    so we conclude thanks to Lemma \ref{lemma: sol of ndKFP are AC}.

    \textbf{Step 2}: on $\R^\infty$ we define, component-wisely, the vector field 
    \begin{equation}
        v_t^{(k)}(\underline{x}):= \begin{cases}
            \int_{\R^d} b_t(x,\mu)\cdot \nabla \varphi_k(x) + \frac{1}{2}a_t(x,\mu):\nablasquare\varphi_k(x) d\mu(x) \quad & \text{ if } \underline{x} = \iota(\mu) \in \mathcal{I}_t
            \\
            0 & \text{ otherwise,}
        \end{cases}
    \end{equation}
    where
    \[\mathcal{I}_t:= \left\{ \underline{x}\in \R^\infty : \underline{x} = \iota(\mu) \in \iota(\PP(\R^d)), \  \int_{\R^d} \frac{|b(t,x,\mu)|}{1+|x|} + \frac{|a(t,x,\mu)|}{1+|x|^2}d\mu(x) <+\infty \right\}.\]
    Note that it is Borel measurable thanks to Corollary \ref{cor: closed and isometry} and well-defined $\widetilde{M}_t$-a.e. for a.e. $t\in(0,T)$ because of the integrability condition \eqref{eq: weighted int cond random measures}. Now, for any $\xi \in C_c^1(0,T)$ and $\tilde{F}\in \operatorname{Cyl}_b^1(\R^\infty)$, i.e. $\tilde{F}(\underline{x}) = \Psi(x_1,\dots,x_k)$ for some $k \in \N$ and $\Psi \in C^1_b(\R^k)$, it holds
    \begin{align*}
        \int_0^T & \xi'(t) \int \tilde{F}(\underline{x})d\widetilde{M}_t(\underline{x}) dt = \int_0^T \xi'(t) \int_{\PP} \Psi(L_{\phi_1(\mu)},\dots, L_{\phi_k}(\mu)) dM_t(\mu)dt 
        \\
        = & - \int_0^T \xi(t) \int \sum_{i=1}^k\partial_i\Psi(L_{\Phi}(\mu)) \int_{\R^d} b_t(x,\mu)\cdot \nabla \phi_i(x) + \frac{1}{2}a_t(x,\mu):\nablasquare\phi_i(x)d\mu(x) dM_t(\mu)dt 
        \\
        = & 
        -\int_0^T \xi(t) \int \sum_{i=1}^k\partial_i\Psi(L_{\Phi}(\mu)) v_t^{(i)}(\iota(\mu)) dM_t(\mu)dt 
        \\
        = & - \int_0^T \xi(t) \int \nabla \tilde{F}(\underline{x})\cdot v_t(\underline{x}) d\widetilde{M}_t(\underline{x})dt,
    \end{align*}
    which means that $\widetilde{M}_t$ solves the continuity equation $\partial_t\widetilde{M}_t + \operatorname{div}_{\R^\infty}(v_t \widetilde{M}_t) = 0$ (see \cite[Section 7]{ambrosio2014well}. Then there exists a measure $\widetilde{\Lambda} \in \PP(C_T(\R^\infty,\tau_w))$, where $\tau_w$ is the element-wise convergence topology over $\R^\infty$ (see \eqref{eq: element-wise distance}), satisfying 
    \begin{itemize}
        \item $(\e_t)_\sharp \widetilde{\Lambda} = \widetilde{M}_t$ for any $t\in [0,T]$;
        \item $\widetilde{\Lambda}$-a.e. $\Tilde{\gamma} \in AC_T(\R^\infty,\tau_w)$, i.e. each component is in $AC_T(\R)$, and it solves
        \[ \partial_t\Tilde{\gamma}^{(i)} = v_t^{(i)}(\Tilde{\gamma}_t) \quad \forall i\in \N, \ \text{for a.e. }t\in [0,T].\]
    \end{itemize}

    \textbf{Step 3}: we prove that $\widetilde{\Lambda}$-a.e. $\Tilde{\gamma}$ is such that $\Tilde{\gamma}(t) \in \iota(\PP(\R^d))$ for all $t\in[0,T]$. To this aim, it is sufficient to prove two things: 
    \begin{itemize}
        \item $\widetilde{\Lambda}$-a.e. $\Tilde{\gamma}$ is such that $\Tilde{\gamma}(t) \in \iota(\PP(\R^d))$ for all $t\in [0,T]\cap \Q$;
        \item $\widetilde{\Lambda}$-a.e. $\Tilde{\gamma}$ is in $AC_T(\R^\infty,D_\infty)$.
    \end{itemize}
    Then, we conclude simply by completeness of $(\R^\infty,D_\infty)$ and the closedness of $\iota(\PP(\R^d))$ in it, given by Corollary \ref{cor: closed and isometry}. The first part is simply obtained noticing that for all $t\in [0,T]\cap \Q$ it holds
    \[\widetilde{\Lambda}\big( \Tilde{\gamma} \ : \ \Tilde{\gamma}(t) \in \iota(\PP(\R^d)) \big) = \widetilde{M}_t(\iota(\PP(\R^d))) = 1.\]
    Regarding the second statement, for $\widetilde{\Lambda}$-a.e. $\Tilde{\gamma}$ and for all $s,t \in [0,T]$ it holds
    \[D_\infty(\Tilde{\gamma}(t),\Tilde{\gamma}(s)) = \sup_{n\in \N} |\Tilde{\gamma}_n(t) - \Tilde{\gamma}_n(s)|\wedge 2 \leq \sup_{n\in \N}\int_s^t |v_r^{(n)}(\Tilde{\gamma}(r))|dr \leq \int_s^t \sup_{n\in \N} |v_r^{(n)}(\Tilde{\gamma}(r))|dr\]
    and it holds
    \[
    \begin{aligned}
    \int & \int_0^T \sup_{n\in \N} |v_r^{(n)}(\Tilde{\gamma}(r))|dr d\widetilde{\Lambda}(\Tilde{\gamma}) =\int_0^T \int \sup_{n\in \N} |v_r^{(n)}(\underline{x})| d\widetilde{M}_r(\underline{x})dr 
    \\
    \leq & \int_0^T \int \sup_{n\in \N} \int_{\R^d} \big|b_r(x,\mu)\cdot \nabla \varphi_n(x)\big| + \frac{1}{2}\big|a_r(x,\mu): \nablasquare\varphi_n(x)\big| d\mu(x)
    \\
    \leq & \int_0^T \int \int \frac{|b_r(x,\mu)|}{1+|x|} + \frac{|a_r(x,\mu)|}{1+|x|^2} d\mu(x)dM_r(\mu)dr<+\infty. 
    \end{aligned}
    \]
    In particular, for $\widetilde{\Lambda}$-a.e. $\Tilde{\gamma}$, the term $ \int_0^T \sup_n |v_r^{(n)}(\Tilde{\gamma}(r))|dr <+\infty$, which concludes the proof of the claim.

    \textbf{Step 4}: the function 
    \[
    \begin{aligned}
    \Theta: C_T(\iota(\PP(\R^d)),D_\infty) \to & C_T(\PP(\R^d)),
    \\
    \Tilde{\gamma}\mapsto &[t\mapsto \iota^{-1}(\Tilde{\gamma}(t))]
    \end{aligned}
    \]
    is well-defined and Borel measurable thanks to Corollary \ref{cor: closed and isometry}. Thus, we define $\Lambda:= \Theta_\sharp  \widetilde{\Lambda} \in \PP(C_T(\PP(\R^d)))$ and verify that it satisfies the requirements. The fact that $(\mathfrak{e}_t)_\sharp \Lambda = M_t$ is straightforward from the definition of $\widetilde{M}_t$ and $(\e_t)_\sharp \widetilde{\Lambda} = \widetilde{M}_t$. Then, it follows  
    \[\int \int_0^T \int \frac{|b_t(x,\mu_t)|}{1+|x|}+\frac{|a_t(x,\mu_t)|}{1+|x|^2} d\mu_t(x) dt d\Lambda(\boldsymbol{\mu}) = \int_0^T \int \int \frac{|b_t(x,\mu)|}{1+|x|}+\frac{|a_t(x,\mu)|}{1+|x|^2} d\mu(x) dM_t(\mu) dt <+\infty.\]
    Regarding the last part, let $\xi \in C_c^1(0,T)$ and $\varphi_k \in \mathcal{S}_{2,w}$. Then, consider a function $\Psi \in C_c^1(\R)$ such that $\Psi(x) = x$ for all $x\in [-1,1]$, so that, defining $\tilde{F}(\underline{x}) = \Psi(x_k) \in \operatorname{Cyl}_c^1 (\R^\infty)$, it holds
    \begin{align*}
        \int \bigg| & \int_0^T \xi'(t) \int \varphi_k(x)d\mu_t(x) dt + \int_0^T \xi(t) \int L_{b_t,a_t}^{\boldsymbol{\mu}}\varphi_k(x) d\mu_t(x) dt\bigg|d\Lambda(\boldsymbol{\mu})
        \\
        = &
        \int \bigg| \int_0^T \xi'(t) \Psi\left(\int \varphi_k(x)d\mu_t(x)\right) dt 
        \\
        & + \int_0^T \xi(t) \Psi'(L_{\varphi_k}(\mu_t)) \int b_t(x,\mu_t)\cdot \nabla\varphi_k(x)+\frac{1}{2}a_t(x,\mu_t):\nablasquare \varphi_k(x) d\mu_t(x) dt\bigg|d\Lambda(\boldsymbol{\mu})
        \\
        = & \int \left| \int_0^T \xi'(t) \tilde{F}(\Tilde{\gamma}(t)) dt + \int_0^T \xi(t) \nabla \tilde{F}(\Tilde{\gamma}(t)) \cdot v_t(\Tilde{\gamma}(t))dt\right|d\widetilde{\Lambda}(\Tilde{\gamma}) = 0.
    \end{align*}
    Selecting $\xi \in \mathcal{A}\subset C_c^1(0,T)$, with $\mathcal{A}$ dense subset in the unit ball of $C_0^1(0,T)$ w.r.t. the norm $\|\cdot\|_{C^1}$, this implies that for $\Lambda$-a.e. $\boldsymbol{\mu} = (\mu_t)_{t\in [0,T]}$, it holds
    \[\int_0^T \xi'(t) \int \varphi(x)d\mu_t(x) dt = - \int_0^T \xi(t) \int L_{b_t,a_t}^{\boldsymbol{\mu}} \varphi(x) d\mu_t(x) dt \quad \forall \xi\in \mathcal{A}, \ \forall \varphi \in \mathcal{S}_{2,w}.\]
    By density of $\mathcal{A}$ in $C_{0}^1(0,T)$ and $\mathcal{S}_{2,w}$ in the unit ball of $C_{0,w}^2(\R^d)$ (that in particular contains the unit ball of $C_c^2(\R^d)$), for $\Lambda$-a.e. $\boldsymbol{\mu}$ it holds $\partial_t \mu_t = (L_{b_t,a_t}^{\boldsymbol{\mu}})^*\mu_t$.
\end{proof}


\subsection{Nested superposition principle: from \texorpdfstring{$\Lambda$}{} to \texorpdfstring{$\mathfrak{L}$}{}}\label{subsec: nestd superposition}
The goal here is to invert the result of Proposition \ref{prop: from L to Lambda nd}, defining a measure $\mathfrak{L}\in \PP(\PP(C_T(\R^d)))$ concentrated over $\MP(b,a)$ given a measure $\Lambda \in \PP(C_T(\PP(\R^d)))$ concentrated over $\KFP(b,a)$, so that $E_\sharp \mathfrak{L} = \Lambda$. The strategy is to use a measurable selection argument to define a map $G:\operatorname{KFP}(b,a) \to \operatorname{MP}(b,a)$ that is a right-inverse for $E$ and that we can use to define $\mathfrak{L}:=G_\sharp \Lambda$. To do so, it is crucial to prove the Borel measurability of the subsets $\KFP(b,a)\subset C_T(\PP(\R^d))$ and $\MP(b,a)\subset \PP(C_T(\R^d))$.

\begin{prop}\label{prop: meas of KFP(b,a)}
    Let $b:[0,T]\times \R^d \times \PP(\R^d) \to \R^d$ and $a:[0,T]\times \R^d \times \PP(\R^d) \to \operatorname{Sym}_+(\R^{d\times d})$ be Borel measurable functions. Then, the subset $\KFP(b,a)\subset C_T(\PP(\R^d))$ is Borel.
\end{prop}

\begin{proof}
    Let $Y:=[0,T]\times \R^d \times \PP(\R^d)$ and define
    \begin{equation}
    \begin{aligned}
        \hat{B}:= \bigg\{ 
            (\hat{\mu}, \hat{\nu},\hat{\eta}) \in \mathcal{M}_+(Y) \times \mathcal{M}(Y;\R^d) \times \mathcal{M}(Y;\R^{d\times d}) \ \text{s.t.} \
            \frac{b}{1+|x|} \in L^{1}(\hat{\mu};\R^d), &  
            \\ \frac{a}{1+|x|^2} \in L^1(\hat{\mu};\R^{d\times d}), \ 
            \hat{\nu} = \frac{b}{1+|x|}\hat{\mu}, \  
            \hat{\eta} = \frac{a}{1+|x|^2} \hat{\mu} &\bigg\}
        \end{aligned}
    \end{equation}
    and
    \begin{equation}\label{eq: hat kfp}
        \begin{aligned}
            &\widehat{\KFP}(b,a):= \bigg\{
            (\hat{\mu}, \hat{\nu},\hat{\eta}) \in \hat{B} \text{ s.t. }
            \forall \xi \in C_c^1(0,T), \ \forall\phi \in C_c^2(\R^d)
            \ \int_Y \xi'(t)\phi(x)d\hat{\mu}(t,x,\mu) =
            \\
            & \ = - \int_Y \xi(t) (1+|x|)\nabla\phi(x)\cdot d\hat{\nu}(t,x,\mu) 
            - \frac{1}{2}\int_Y \xi(t)(1+|x|^2)\nablasquare\phi(x):d\hat{\eta}(t,x,\mu) 
            \bigg\}.
        \end{aligned}
    \end{equation}
    We endow $\mathcal{M}_+(Y)$ and $\mathcal{M}(Y;\R^n)$ with the narrow topology, so that $\mathcal{M}_+(Y) \times \mathcal{M}(Y;\R^d) \times \mathcal{M}(Y;\R^{d\times d})$ is endowed with the product topology. Then, we claim that $\widehat{\operatorname{KFP}}(b,a)$ is a Borel subset of the product: indeed the integral equality is a closed condition, since all the integrals involved in \eqref{eq: hat kfp} are continuous functions, and $\hat{B}$ is Borel thanks to Lemma \ref{lemma: meas results}.
    Now, define the maps 
    \begin{equation}
    \begin{aligned}
        & \pi^1: \mathcal{M}_+(Y) \times \mathcal{M}(Y;\R^d) \times \mathcal{M}(Y;\R^{d\times d}) \to \mathcal{M}_+(Y), \quad \pi^1(\hat{\mu},\hat{\nu},\hat{\eta})  = \hat{\mu}
        \\
        & \kappa:C_T(\PP(\R^d)) \to \mathcal{M}_+(Y), \quad \kappa(\boldsymbol{\mu}) = dt \otimes (\mu_t\otimes \delta_{\mu_t}) = \int_0^T \delta_t\otimes \mu_t \otimes \delta_{\mu_t} dt
    \end{aligned}
    \end{equation}
    The projection map $\pi^1|_{\hat{B}}$ is continuous and injective, thus it maps Borel sets in Borel sets (see e.g. \cite{bogachev2007measure} for a detailed description or \cite[Appendix A]{pinzisavare2025} for a quick overview). It can be easily seen that the function $\kappa$ is continuous and injective. 
    Then, we conclude proving that $\KFP(b,a) = \kappa^{-1}\big(\pi^1(\hat{\KFP}(b,a))\big)$:
    \begin{itemize}
        \item let $(\mu_t)_{t\in [0,T]}\in \KFP(b,a)$, then $(\hat{\mu},\hat{\nu},\hat{\eta}):=\big(\kappa(\boldsymbol{\mu}), \frac{b}{1+|x|}\kappa(\boldsymbol{\mu}), \frac{a}{1+|x|^2}\kappa(\boldsymbol{\mu}) \big) \in \widehat{\KFP}(b,a)$, indeed clearly the conditions on the densities are verified and 
        \[\int_Y \frac{|b(t,x,\mu)|}{1+|x|} + \frac{|a(t,x,\mu)|}{1+|x|^2} d\big(\kappa(\boldsymbol{\mu})\big)(t,x,\mu) = \int_0^T \int_{\R^d}\frac{|b(t,x,\mu_t)|}{1+|x|} + \frac{|a(t,x,\mu_t)|}{1+|x|^2}d\mu_t(x) dt<+\infty,\]
        \[
        \begin{aligned}
        0 = & \int_0^T \xi'(t)\int \phi(x)d\mu_t(x)dt+ \int_0^T \xi(t)\int \nabla\phi(x)\cdot b(t,x,\mu_t) d\mu_t(x) dt 
        \\
        & + \frac{1}{2}\int_0^T \xi(t)\int \nablasquare\phi(x): a(t,x,\mu_t)d\mu_t(x)dt
        \\
        = & \int_0^T \int_{\PP}\int_{\R^d} \xi'(t)\phi(x) d\mu(x) d\delta_{\mu_t}(\mu) dt + \int_0^T \int_{\PP}\int_{\R^d}\xi(t) \nabla\phi(x) \cdot b(t,x,\mu) d\mu(x) d\delta_{\mu_t}(\mu) dt
        \\
        & + \frac{1}{2} \int_0^T \int_{\PP}\int_{\R^d} \xi(t)\nablasquare\phi(x) : a(t,x,\mu) d\mu_t(x) d\delta_{\mu_t}(\mu) dt
        \\
        = & \int_Y \xi'(t)\phi(x) d\big(\kappa(\boldsymbol{\mu})\big)(t,x,\mu) + \int_Y \xi(t) \nabla\phi(x) \cdot b(t,x,\mu) d\big(\kappa(\boldsymbol{\mu})\big)(t,x,\mu) 
        \\
        & + \frac{1}{2}\int_Y \xi(t)\nablasquare\phi(x):a(t,x,\mu)d\big(\kappa(\boldsymbol{\mu})\big)(t,x,\mu)
        \\
        = & 
        \int_Y \xi'(t)\phi(x) d\hat{\mu}(t,x,\mu) + \int_Y \xi(t) (1+|x|)\nabla\phi(x) \cdot d\hat{\nu}(t,x,\mu) 
        \\
        & + \frac{1}{2}\int_Y \xi(t)(1+|x|^2)\nablasquare\phi(x):d\hat{\eta}(t,x,\mu),
        \end{aligned}
        \]
        for all $\xi\in C_c^1(0,T)$ and $\phi \in C_c^2(\R^d)$;
        \item vice versa, let $\boldsymbol{\mu} \in \kappa^{-1}\big(\pi^1(\widehat{\KFP}(b,a))\big)$, which means that $\big(\kappa(\boldsymbol{\mu}), \frac{b}{1+|x|}\kappa(\boldsymbol{\mu}), \frac{a}{1+|x|^2}\kappa(\boldsymbol{\mu}) \big) \in \widehat{\KFP}(b,a)$, and we conclude that $\boldsymbol{\mu}\in \KFP(b,a)$ thanks to the same computations made in the previous case.
    \end{itemize}
    \vspace{-0.5cm}
\end{proof}

\vspace{0.1cm}

\begin{prop}\label{prop: meas of MP}
     Let $b:[0,T]\times \R^d \times \PP(\R^d) \to \R^d$ and $a:[0,T]\times \R^d \times \PP(\R^d) \to \operatorname{Sym}_+(\R^{d\times d})$ be Borel measurable functions. Then, the subset $\MP(b,a)\subset \PP(C_T(\R^d))$ is Borel.
\end{prop}

\begin{proof}
    \textbf{Step 1}: let us start by fixing the notation. Let $\lambda \in \MP(b,a)$, by \eqref{eq: MP(b,a)} and Definition \ref{def: MP solution fin dim}, for all $\xi\in C_c^1(0,T)$ and $\phi\in C_c^2(\R^d)$ it holds
    \[t\mapsto X_t^{\xi\phi,\lambda}(\gamma) = X_t^{\xi\phi}(\gamma,\lambda) := \xi(t)\phi(\gamma_t) - \int_0^t \xi'(r)\phi(\gamma_r) + \xi(r)L_{b_r,a_r}^\lambda \phi(\gamma_r)dr\]
    is a martingale in the filtered space $\big(C_T(\R^d), (\mathcal{F}_t)_{t\in [0,T]},\mathcal{F}_T, \lambda\big)$. The integrability conditions ensure that $X_t^{\xi\phi}\in L^1(\lambda)$, so that the martingale condition can be rewritten as follows: 
    \begin{equation}\label{eq: mart property proof meas}
    \begin{gathered}
    \text{$\forall\xi\in C_c^1(0,T)$, 
     $\ \forall\phi\in C_c^2(\R^d)$, $ \ \forall \, 0\leq s < t \leq T$, $ \ \forall A\in \mathcal{F}_s$}
    \\
    \int_{C_T(\R^d)} \mathds{1}_A(\gamma) \big(X_{t}^{\xi\phi,\lambda}(\gamma) - X_{s}^{\xi\phi,\lambda}(\gamma)\big) d\lambda(\gamma) = 0.
    \end{gathered}
    \end{equation}
    To prove measurability of $\MP(b,a)$, it is first necessary to show that such a condition can be asked to hold for a countable number of $\xi,\phi,s,t,A$. Before proceeding, it will be useful the following (uniform in time) estimate, that follows from \eqref{eq: estimate for the operator}:
    \begin{equation}\label{eq: estimate for martingale}
        \begin{aligned}
            | & X_t^{\xi\phi,\lambda}(\gamma)| \leq (T+1)\|\xi\|_{C^1} \|\phi\|_\infty 
            \\
            & + \|\xi\|_\infty \int_0^T |b(r,\gamma_r,(\e_r)_\sharp \lambda)||\nabla\phi| + |a(r,\gamma_r,(\e_r)_\sharp \lambda)| |\nablasquare \phi| dr 
            \\
            \leq & (T+3)\|\xi\|_{C^1}\|\phi\|_{C^2_{0,w}}\left( 1 + \left\|\frac{b^{\lambda}(\cdot,\gamma)}{1+|\gamma(\cdot)|}\right\|_{L^1(0,T)} + \left\|\frac{a^\lambda(\cdot,\gamma)}{1+|\gamma(\cdot)|^2}\right\|_{L^1(0,T)}\right),
        \end{aligned}
    \end{equation}
    where $b^\lambda(t,\gamma):= b(t,\gamma_t,(\e_t)_\sharp \lambda)$ and similarly for $a^\lambda$. In particular, for every $\xi\in C_c^1(0,T)$ and $\phi\in C_c^2(\R^d)$, the right hand side is finite for $\lambda$-a.e. $\gamma$ because $\lambda \in \operatorname{MP}(b,a)$. 

    \textbf{Step 2}: it is sufficient to check \eqref{eq: mart property proof meas} for a countable amount of $\xi$ and $\phi$. Indeed, consider $\mathcal{D}_T\subset C_0^1(0,T)$ and $\mathcal{D}_{2,w}\subset C_{0,w}^2(\R^d)$ dense and countable subsets. Notice that, by linearity of the integral and of the operator $L_{b_t,a_t}^\lambda$, for any $\xi_1,\xi_2 \in C_c^1(0,T)$ and $\phi_1,\phi_2 \in C_c^2(\R^d)$ it holds
    \[X^{\xi_1\phi_1 +\xi_2\phi_2,\lambda}_t(\gamma) = X^{\xi_1\phi_1,\lambda}_t(\gamma)+\X^{\xi_2\phi_2,\lambda}_t (\gamma) \quad \forall \lambda \in \PP(C_T(\R^d)), \ \gamma \in C_T(\R^d).\]
    Assume that \eqref{eq: mart property proof meas} holds for all $\xi \in \mathcal{D}_T$ and $\phi\in \mathcal{D}_{2,w}$. Then, for any $\xi \in C_c^1(0,T)$ and $\phi\in C_c^2(\R^d)$, consider $\xi_n \in \mathcal{D}_T$ and $\phi_n\in \mathcal{D}_{2,w}$ such that $\xi_n \to \xi$ and $\phi_n \to \phi$ in their respective spaces, so that for all $0\leq s<t \leq T$ and $A \in \mathcal{F}_s$, it holds 
    \begin{align*}
        \bigg|\int \mathds{1}_A(\gamma) & \big(  X_{t}^{\xi\phi,\lambda}(\gamma) - X_{s}^{\xi\phi,\lambda}(\gamma)\big) d\lambda(\gamma) \bigg| 
        =  \left|\int \mathds{1}_A(\gamma) \big(X_{t}^{\xi\phi-\xi_n\phi_n,\lambda}(\gamma) - X_{s}^{\xi\phi -\xi_n\phi_n,\lambda}(\gamma)\big) d\lambda(\gamma) \right|\
        \\
        & \leq \int |X_t^{\xi(\phi-\phi_n) + (\xi-\xi_n)\phi_n}(\gamma)|+|X_s^{\xi(\phi-\phi_n) + (\xi-\xi_n)\phi_n}(\gamma)|d\lambda(\gamma)
        \\
        & \leq (T+3)\left(\Big(\sup_{n\in \N}\|\phi_n\|_{C_{0,w}^2}\Big)\|\xi_n- \xi\|_{C^1}+ \|\xi\|_{C^1}\|\phi_n - \phi\|_{C^2_{0,w}}\right)\cdot
        \\
        & \ \ \cdot\bigg( 1 + \int \int_0^T\frac{|b(t,\gamma_t,(\e_t)_\sharp \lambda)|}{1+|\gamma_t|} + \frac{|a(t,\gamma_t,(\e_t)_\sharp \lambda)|}{1+|\gamma_t|^2} dt d\lambda(\gamma)\bigg) \to 0.
    \end{align*}

    \textbf{Step 3}: it is sufficient to check \eqref{eq: mart property proof meas} for $s<t$ with $s, t \in \mathbb{Q}\cap[0,T]$. Indeed, consider any $0\leq s<t \leq T$ and assume that \eqref{eq: mart property proof meas} holds for two sequences of times $s_n <t_n$ satisfying $s_n \searrow s$, $t_n \to t$. Then, consider any $\xi\in\mathcal{D}_T$, $\phi \in \mathcal{D}_{2,w}$ and $A \in \mathcal{F}_s$, so that $A$ is also $\mathcal{F}_{s_n}$-measurable for any $n\in \N$. By continuity of $\xi$ and $\phi$, and by $\int_0^T |L_{b_r,a_r}^\lambda \phi(\gamma_r)|dr<+\infty$ for $\lambda$-a.e. $\gamma$, it holds that $X_{s_n}^{\xi\phi,\lambda} \to X_s^{\xi\phi,\lambda}$ and $X_{t_n}^{\xi\phi,\lambda} \to X_{t}^{\xi\phi,\lambda}$ for $\lambda$-a.e. $\gamma$, dominated by the right-hand side of \eqref{eq: estimate for martingale}, which is in $L^1(\lambda)$. Finally, we can apply dominated convergence theorem to obtain 
    \[0 = \int_{C_{T}(\R^d)} \mathds{1}_A(\gamma) \big(X_{t_n}^{\xi\phi,\lambda}(\gamma) - X_{s_n}^{\xi\phi,\lambda}(\gamma)\big) d\lambda(\gamma) \to \int_{C_{T}(\R^d)}\mathds{1}_A(\gamma) \big(X_{t}^{\xi\phi,\lambda}(\gamma) - X_{s}^{\xi\phi,\lambda}(\gamma)\big) d\lambda(\gamma).\]

    \textbf{Step 4}: it is sufficient to check \eqref{eq: mart property proof meas} for a countable amount of $A \in \mathcal{F}_s$. Indeed, fix $\xi \in \mathcal{D}_T$, $\phi \in \mathcal{D}_{2,w}$ and $0\leq s < t \leq T$ with $s,t\in \Q$. Then, thanks to Lemma \ref{lemma: equiv of sigma-algebras}, there exists a countable algebra $\mathcal{V}_s \subset \mathcal{F}_s$ that generates the whole $\sigma$-algebra. We check that it is sufficient to check \eqref{eq: mart property proof meas} for $A\in \mathcal{V}_s$. Indeed the class
    \[\left\{ A \in \mathcal{F}_s : \int_{C_T(\R^d)} \mathds{1}_A \big( X_{t}^{\xi\phi,\lambda}(\gamma) - X_{s}^{\xi\phi,\lambda}(\gamma) \big) d\lambda(\gamma) = 0 \right\}\]
    is a Dynkin system and it contains the algebra $\mathcal{V}_s$. Then, by Dynkin's lemma and $L^1$-integrability of $X_{t}^{\xi\phi,\lambda} - X_{s}^{\xi\phi,\lambda}$, it contains its generated $\sigma$-algebra $\mathcal{F}_s$, concluding this step.

    \textbf{Step 5}: to recap, we can write $\MP(b,a)$ as 
    \begin{equation}\label{eq: intersection mp}
        \MP(b,a) = \bigcap_{\substack{\xi\in \mathcal{D}_T, \\ \phi\in \mathcal{D}_{2,w}}} \bigcap_{\substack{s,t \in [0,T]\cap \mathbb{Q}, \\ s<t}} \bigcap_{A\in \mathcal{V}_s} \MP(b,a; \xi,\phi,s,t,A), 
    \end{equation}
    where 
    \begin{equation}
    \begin{aligned}
        \MP(b,a; \xi,\phi,s,t,A) := \bigg\{\lambda \in \PP(C_T(\R^d)) \ : \hspace{2cm} & \\
        \int\int_0^T \frac{|b_t(\gamma_t,(\e_t)_\sharp \lambda) |}{1+|\gamma_t|} + \frac{|a_t(\gamma_t,(\e_t)_\sharp \lambda)|}{1+|\gamma_t|^2} dtd\lambda(\gamma) <+\infty, & 
        \\
        \int_{C_T(\R^d)} \mathds{1}_A(\gamma) \big( X_t^{\xi\phi,\lambda}(\gamma) - X_s^{\xi\phi,\lambda}(\gamma) \big) d\lambda(\gamma) = 0 &
        \bigg\}.
    \end{aligned}
    \end{equation}
    So, we are left with the proof of Borel measurability of the sets $\MP(b,a;\xi,\phi,s,t,A)$. Define the Polish space $Z:= C_T(\R^d) \times \PP(C_T(\R^d))$ endowed with the product topology, and the natural injection from $\PP(C_T(\R^d))$ to $\PP(Z)$ as
    \begin{equation}
        \mathfrak{K}: \PP(C_T(\R^d)) \to \PP(Z),\quad \mathfrak{K}(\lambda) := \lambda \otimes \delta_{\lambda}.
    \end{equation}
    Then, for any $\xi$, $\phi$, $s$, $t$ and $A$ as in \eqref{eq: intersection mp} define 
    \[
    \begin{aligned}
    \widehat{\MP}:= \bigg\{\hat{\lambda} \in \mathcal{M}_+(Z) \ : \ & \int_Z \int_0^T\frac{|b(r,\gamma_r,(\e_r)_\sharp \lambda) |}{1+|\gamma_r|} + \frac{|a(r,\gamma_r,(\e_r)_\sharp \lambda)|}{1+|\gamma_r|^2} dr \ d\hat{\lambda}(\gamma,\lambda) <+\infty, 
    \\
    & \int_{Z} \mathds{1}_A(\gamma) \big( X_t^{\xi\phi,\lambda}(\gamma) - X_s^{\xi\phi,\lambda}(\gamma) \big) d\hat{\lambda}(\gamma,\lambda) = 0\bigg\}.
    \end{aligned}
    \]
    Since, for all $r\in [0,T]$, the map $Z \ni (\gamma,\lambda)\mapsto (\gamma_r,(\e_r)_\sharp \lambda)$ is continuous (and in particular measurable) for all $r\in [0,T]$, then also the map
    \[(\gamma,\lambda)\mapsto \int_0^T \frac{|b(r,\gamma_r,(\e_r)_\sharp \lambda) |}{1+|\gamma_r|} + \frac{|a(r,\gamma_r,(\e_r)_\sharp \lambda)|}{1+|\gamma_r|^2} dr\]
    is measurable (see, e.g., \cite[§3.4]{bogachev2007measure}). Moreover, looking at the definition of $X_t^{\xi\phi,\lambda}(\gamma)$, it is not hard to realize that also
    \[(\gamma,\lambda) \mapsto X_t^{\xi\phi,\lambda}(\gamma)\]
    is measurable for any $t\in [0,T]$. Then, thanks to Lemma \ref{lemma: meas results}, the set $\widehat{\MP}$ is Borel. Then, we conclude noticing that $\MP(b,a;\xi,\phi,s,t,A) = \mathfrak{K}^{-1}(\widehat{\MP})$.
\end{proof}

We are ready to prove the main theorem of this section, that we call \textit{nested superposition principle for SDE}.

\begin{teorema}[Nested superposition principle for SDE]\label{thm: from Lambda to mathfrak L}
     Let $b:[0,T]\times \R^d \times \PP(\R^d) \to \R^d$ and $a:[0,T]\times \R^d \times \PP(\R^d) \to \operatorname{Sym}_+(\R^{d\times d})$ be Borel measurable. Then there exists a Souslin-Borel measurable map $G_{b,a}: \operatorname{KFP}(b,a) \to \PP(C_T(\R^d))$ satisfying $\operatorname{Im}(G_{b,a})\subset \MP(b,a)$ and $E\circ G_{b,a}(\boldsymbol{\mu}) = \boldsymbol{\mu}$ for all $\boldsymbol{\mu} \in \KFP(b,a)$, i.e. $G_{b,a}$ is a right-inverse for $E|_{\MP(b,a)}$.
     \\
     In particular, if $\Lambda\in \PP(C_T(\PP(\R^d)))$ is concentrated over $\KFP(b,a)$, then $\mathfrak{L}:=(G_{b,a})_\sharp \Lambda \in \PP(\PP(C_T(\R^d)))$ is well-defined, it is concentrated over $\MP(b,a)$ and it satisfies $E_\sharp \mathfrak{L} = \Lambda$.
\end{teorema}

\begin{proof}
    The restriction map 
    \[E|_{\MP(b,a)} : \MP(b,a) \to \KFP(b,a)\]
    is well-defined, thanks to Proposition \ref{prop: from L to Lambda nd}. Moreover, because of the finite dimensional superposition for SDE, i.e. Theorem \ref{thm: fin dim superposition}, it is surjective. Then, thanks to Propositions \ref{prop: meas of KFP(b,a)} and \ref{prop: meas of MP}, we can apply \cite[Theorem 6.9.1]{bogachev2007measure} to obtain a Souslin-Borel measurable map $G_{b,a}:\KFP(b,a)\to \PP(C_T(\R^d))$ satisfying the requirements. Then, thanks to the universal measurability of Souslin sets (see \cite[Chapter 6]{bogachev2007measure}), the measure $\mathfrak{L}:= (G_{b,a})_\sharp \Lambda$ is a well-defined Borel measure, and by the properties of $G_{b,a}$ it satisfies the requirements.
\end{proof}


Finally, putting all the results of this section together, we have a proof for the nested stochastic superposition principle, Theorem \ref{main theorem}.

\begin{proof}[Proof of Theorem \ref{main theorem}]
    The existence of $\Lambda$ and property (1) come from Theorem \ref{thm: superposition from M to Lambda}. Then, the existence of $\mathfrak{L}$, together with the properties (2) and (3), is a consequence of Theorem \ref{thm: from Lambda to mathfrak L}. On the other hand, (i) and (ii) follow, respectively, from Proposition \ref{prop: from L to Lambda nd} and Proposition \ref{prop: from Lambda to M nd}.
\end{proof}
\section{Uniqueness scheme}\label{sec: uniqueness}
In this section, we show how uniqueness can be transferred between the main objects of Theorem \ref{main theorem}, $\boldsymbol{M}\in C_T(\PP(\PP(\R^d)))$, $\Lambda \in \PP(C_T(\PP(\R^d)))$ and $\mathfrak{L}\in \PP(\PP(C_T(\R^d)))$.

A first result in this direction is a consequence of the nested superposition principle for SDEs.

\begin{lemma}\label{lemma: transfer of uniqueness}
	Let $b:[0,T]\times \R^d \times \PP(\R^d) \to \R^d$ and $a:[0,T]\times \R^d \times \PP(\R^d) \to \operatorname{Sym}_+(\R^{d\times d})$ Borel measurable maps. Let $\overline{M} \in \PP(\PP(\R^d))$. The following are equivalent:
	\begin{enumerate}
		\item there exists at most one curve $\boldsymbol{M} = (M_t)_{t\in[0,T]} \in C_T(\PP(\PP(\R^d)))$ satisfying \eqref{eq: equation main theorem} and $M_0 = \overline{M}$;
		\item there exists at most one $\Lambda \in \PP(C_T(\PP(\R^d)))$ concentrated on $\operatorname{KFP}(b,a)$ and satisfying \eqref{eq: integr of Lambda} and $(\mathfrak{e}_0)_\sharp \Lambda = \overline{M}$;
		\item given $\mathfrak{L}^1,\mathfrak{L}^2 \in \PP(\PP(C_T(\R^d)))$ concentrated on $\operatorname{MP}(b,a)$ and satisfying \eqref{eq: integr of mathfrak L} and $(E_0)_\sharp \mathfrak{L}^1 = (E_0)_\sharp \mathfrak{L}^2 = \overline{M}$, then $E_\sharp \mathfrak{L}^1 = E_\sharp \mathfrak{L}^2$ as elements of $\PP(C_T(\PP(\R^d)))$.
	\end{enumerate}
	In particular, if one of the above points holds together with existence, then existence for the other points holds, and for $\overline{M}$-a.e. $\overline{\mu}\in \PP(\R^d)$ the following hold:
	\begin{itemize}
		\item[(a)] there exists a unique curve of measures $\boldsymbol{\mu} = (\mu_t)_{t\in[0,T]}\in \operatorname{KFP}(b,a)$ with $\mu_0 = \overline{\mu}$;
		\item[(b)] there exists $\lambda \in \operatorname{MP}(b,a)$ with $(\e_0)_\sharp \lambda = \overline{\mu}$ and, given $\lambda^1,\lambda^2 \in \operatorname{MP}(b,a)$ satisfying $(\e_0)_\sharp \lambda^1 = (\e_0)_\sharp \lambda^2 = \overline{\mu}$, then $(\e_t)_\sharp \lambda^1 = (\e_t)_\sharp \lambda^2$ for all $t\in [0,T]$.
	\end{itemize}
\end{lemma}

\begin{proof}
	(2)$\implies$(3): thanks to Proposition \ref{prop: from L to Lambda nd}, $\Lambda^1 :=E_\sharp \mathfrak{L}^1$ and $\Lambda^2:=E_{\sharp }\mathfrak{L}^2$ are concentrated on $\operatorname{KFP}(b,a)$, and because of \eqref{eq: integr of mathfrak L}, they satisfy \eqref{eq: integr of Lambda} as well. Then $\Lambda^1 = \Lambda^2$.
	
	(3)$\implies$(2): assume that $\Lambda^1,\Lambda^2 \in \PP(C_T(\PP(\R^d)))$ satisfy the conditions in (2). Thanks to Theorem \ref{thm: from Lambda to mathfrak L}, there exist $\mathfrak{L}^1,\mathfrak{L}^2 \in \PP(\PP(C_T(\R^d)))$ satisfying the conditions in (3), and we conclude by observing that $\Lambda^1 = E_\sharp \mathfrak{L}^1 = E_\sharp \mathfrak{L}^2 = \Lambda^2$.
	
	(2)$\implies$(1): let $\boldsymbol{M}^1, \boldsymbol{M}^2 \in C_T(\PP(\PP(\R^d)))$ satisfy the conditions in (1). We can lift both to $\Lambda^1,\Lambda^2 \in \PP(C_T(\PP(\R^d)))$ thanks to Theorem \ref{thm: superposition from M to Lambda}. It is easy to see that both the liftings satisfy the conditions in (2), so that they coincide, from which $M_t^1 = (\mathfrak{e}_t)_\sharp \Lambda^1 = (\mathfrak{e}_t)_\sharp \Lambda^2 = M_t^2$ follows for all $t\in[0,T]$.
	
	(1)$\implies$(2): if no such $\Lambda$ exists, then there is nothing to prove. Then, assume $\Lambda$ satisfies the conditions in (2). Disintegrating it with respect to $\mathfrak{e}_0$, we obtain 
	\[\Lambda = \int_{\PP(\R^d)} \Lambda_{\overline{\mu}}d\overline{M}(\overline{\mu}),\]
	where, for $\overline{M}$-a.e. $\overline{\mu}$, say $\overline{\mu}\in \mathcal{N}^c\in \mathcal{B}(\PP(\R^d))$ with $\overline{M}(\mathcal{N}) = 0$, $\Lambda_{\overline{\mu}} \in \PP(C_T(\PP(\R^d)))$ is concentrated on $\operatorname{KFP}(b,a)\cap\{\boldsymbol{\mu}:\mu_0=\overline{\mu}\}$, which in particular is non-empty. We are done if we show that, for $\overline{M}$-a.e. $\overline{\mu}\in \mathcal{N}^c$, such a set is a singleton. Let $$\mathcal{A}:= \left\{(\boldsymbol{\mu}^1,\boldsymbol{\mu}^2) \in \big(C_T(\PP(\R^d))\big)^2 : \boldsymbol{\mu}^1 \neq \boldsymbol{\mu}^2, \ \boldsymbol{\mu}^1(0) = \boldsymbol{\mu}^2(0) \in \mathcal{N}^c, \ \boldsymbol{\mu}^1,\boldsymbol{\mu}^2 \in \KFP(b,a)\right\}.$$
    It is a Borel set, and the map $\mathfrak{e}_0^1: \mathcal{A} \to \mathcal{N}^c$, $\mathfrak{e}_0^1(\boldsymbol{\mu}^1,\boldsymbol{\mu}^2) = \boldsymbol{\mu}^1(0)$ is well-defined. Since it is continuous, its image $\operatorname{Im}(\mathfrak{e}_0^1) \subset \mathcal{N}^c$ is a Souslin subset of $\PP(\R^d)$. By contradiction, assume $\overline{M}(\operatorname{Im}(\mathfrak{e}_0^1)) > 0$. There exists a Souslin/Borel measurable map (see \cite[Theorem 6.9.1]{bogachev2007measure}) $G = (G^1,G^2): \operatorname{Im}(\mathfrak{e}_0^1) \to \mathcal{A}$ that is a right-inverse of $\mathfrak{e}_0^1$. We impose two further restrictions: for all $n\in \N$ consider the $\overline{M}$-measurable set
    \begin{equation*}
        \mathcal{G}_n:= \left\{ \overline{\mu} \in \operatorname{Im}(\mathfrak{e}_0^1) : i=1,2, \text{ }(\mu_t) = G^i(\overline{\mu}) \Rightarrow \int_0^T \hspace{-0.2cm} \int_{\R^d} \frac{|b_t(x,\mu_t)|}{1+|x|} + \frac{|a_t(x,\mu_t)|}{1+|x|^2} d\mu_t(x) dt \leq n \right\}.
    \end{equation*}
    Since $ G^i(\overline{\mu}) \in \KFP(b,a)$, we have $\cup_{n\in \N} \mathcal{G}_n = \operatorname{Im}(\mathfrak{e}_0^1)$, thus there exists $\bar{n}\in \N$ such that $\overline{M}(\mathcal{G}_{\bar{n}}) >0$.
    
    Moreover, for all $\overline{\mu} \in \mathcal{G}_{\bar{n}}$ there exists $t\in(0,T)\cap \mathds{Q}$ such that $\mathfrak{e}_t(G^1(\overline{\mu})) \neq \mathfrak{e}_t(G^2(\overline{\mu}))$. In particular, there exists $\bar{t}\in(0,T)\cap \mathds{Q}$ such that $\overline{M}(\mathcal{G}_{\bar{n}} \cap \{\mathfrak{e}_{\bar{t}}\circ G^1 \neq \mathfrak{e}_{\bar{t}} \circ G^2\})>0$. By separability of $\PP(\R^d)$, we can find a countable open cover $\{\mathcal{I}^1_{j} \times \mathcal{I}^2_j : j\in \N \}$ for the off-diagonal set $\PP(\R^d)\times \PP(\R^d) \setminus\{(\mu,\mu) : \mu\in  \PP(\R^d)\}$. In particular, $\mathcal{I}^1_j \cap \mathcal{I}^2_j = \emptyset$ for all $j \in \N$. Then, there exists $\bar{j}\in \N$ such that $\overline{M}( \mathcal{H}(\bar{n},\bar{t}, \bar{j} ) ) >0$, where
    \[\mathcal{H} (\bar{n},\bar{t}, \bar{j}) := \mathcal{G}_{\bar{n}} \cap (\mathfrak{e}_{\bar{t}}\circ G^1)^{-1}(\mathcal{I}^1_{\bar{j}}) \cap (\mathfrak{e}_{\bar{t}} \circ G^2)^{-1}(\mathcal{I}^2_{\bar{j}})\]
    is an $\overline{M}$-measurable set.
    For $i=1,2$ define 
    \begin{align*}
        M_t^i = N_t^{i} + \tilde{N}_t :=  \int_{\mathcal{H} (\bar{n},\bar{t}, \bar{j}) } \delta_{\mathfrak{e}_t(G^i(\overline{\mu}))} d\overline{M}(\overline{\mu}) + \int_{\PP(\R^d) \setminus \mathcal{H} (\bar{n},\bar{t}, \bar{j}) } (\mathfrak{e}_t)_\sharp \Lambda_{\overline{\mu}} \, d\overline{M} (\overline{\mu}).
    \end{align*}
    Now, $M_0^i = \overline{M}$ for $i=1,2$ and the two curves satisfy (1): indeed $M_t^i = (\mathfrak{e}_t)_\sharp \Lambda^i$, where 
    \[\Lambda_i := \int_{\mathcal{H} (\bar{n},\bar{t}, \bar{j}) } \delta_{G^i(\overline{\mu})} d\overline{M}(\overline{\mu}) + \int_{\PP(\R^d) \setminus \mathcal{H} (\bar{n},\bar{t}, \bar{j}) }  \Lambda_{\overline{\mu}} \, d\overline{M} (\overline{\mu}),\]
    which satisfies $(\mathfrak{e}_0)_\sharp \Lambda^i = \overline{M}$, \eqref{eq: integr of Lambda} and is concentrated on $\KFP(b,a)$ by construction, thus we can apply Proposition \ref{prop: from Lambda to M nd}. Finally, we reach a contradiction by observing that, by construction $ M_{\bar{t}}^1 (\mathcal{I}_{\bar{j}}^1) = \overline{M}(\mathcal{H} (\bar{n},\bar{t}, \bar{j})) + \tilde{N}_{\bar{t}}(\mathcal{I}_{\bar{j}}^1) > \tilde{N}_{\bar{t}}(\mathcal{I}_{\bar{j}}^1) = M_{\bar{t}}^2 (\mathcal{I}_{\bar{j}}^1)$, so that $ M_{\bar{t}}^1 \neq M_{\bar{t}}^2$.

	We conclude noticing that (a) and (b) are byproducts of the previous argument. \qedhere
\end{proof}

As already observed in \cite[pp. 11]{trevisan2016well}, the presence of the diffusion term does not allow us to conclude uniqueness of martingale solutions only assuming the solution of the associated Kolmogorov-Fokker-Planck equation is uniquely determined by the starting measure. Here, also the non-local nature of the problem is an obstacle for proving the uniqueness of $\mathfrak{L} \in \PP(\PP(C_T(\R^d)))$ assuming uniqueness of $\Lambda \in \PP(C_T(\PP(\R^d)))$. In particular, we need to adapt \cite[Proposition 2.6]{trevisan2016well} to the following.

\begin{lemma}\label{lemma: restriction}
    Let $b:[0,T]\times \R^d \times \PP(\R^d) \to \R^d$ and $a:[0,T]\times \R^d \times \PP(\R^d) \to \operatorname{Sym}_+(\R^{d\times d})$ be Borel measurable maps. Let $0 <s <T$ and $\lambda \in \operatorname{MP}(b,a)$. Then:
    \begin{itemize}
        \item the measure $(|_{[s,T]})_\sharp  \lambda \in \PP(C([s,T],\R^d))$ is a martingale solution associated to $L_{b_t,a_t}^{\lambda}$ in the interval $[s,T]$;
        \item let $\rho:C([0,T],\R^d)\to [0,+\infty)$ be a bounded probability density with respect to $\lambda$ and assume it is $\mathcal{F}_s$-measurable. Then $(\rho\lambda)|_{[s,T]}:=(|_{[s,T]})_\sharp (\rho\lambda) \in \PP(C([s,T],\R^d))$ is a martingale solution associated to the operator $L_{b_t,a_t}^{\lambda}$ in the interval $[s,T]$.
    \end{itemize}
\end{lemma}

\begin{proof}
    It is an immediate consequence of \cite[Proposition 2.6]{trevisan2016well}.
\end{proof}

In view of the previous lemma, we make the following assumption on the coefficients $b$ and $a$. We adopt the obvious notation $\operatorname{KFP}_{[s,T]}(b,a)$ for the curves of probability measures $(\mu_t)_{t\in[s,T]}$ that solves $\partial_t \mu_t = (L_{b_t,a_t}^{\boldsymbol{\mu}})^*\mu_t$ in $[s,T]$ and satisfy the integrability condition in \eqref{eq: KFP(b,a)} integrating between $s$ and $T$. The similar notation $\operatorname{MP}_{[s,T]}(b,a)$ is used for the set of martingale solutions.

\begin{ass}\label{assumption uniqueness}
    For all $s\in [0,T]$, for all $\boldsymbol{\mu} \in \operatorname{KFP}_{[s,T]}(b,a)$, and for all $\tilde{\mu}\ll\mu_s$ such that $\tilde{\mu} \in \PP(\R^d)$ and $\frac{d\tilde{\mu}}{d\mu_s} \in L^\infty(\mu_s)$, there exists at most one $\tilde{\boldsymbol{\mu}} \in C([s,T],\PP(\R^d))$ such that $\tilde{\mu}_s=\tilde{\mu}$, $\tilde{\mu}_t\ll\mu_t$ for all $t\in[s,T]$, $\frac{d\tilde{\mu}_t}{d\mu_t} \in L^\infty(\mu_t)$ and $\partial_t\tilde{\mu}_t = (L_{b_t,a_t}^{\boldsymbol{\mu}})^*\tilde{\mu}_t$.
\end{ass}

We are asking for uniqueness for the linearized Kolmogorov-Fokker-Planck problems associated to $L_{b_t,a_t}^{\boldsymbol{\mu}}$, fixing $\boldsymbol{\mu}\in C_T(\PP(\R^d))$. It is quite a natural strategy to study the uniqueness of the linearized version of the KFP equation to then obtain uniqueness of the non-linear one (see e.g. \cite{barbu2021uniqueness}). The following shows a case in which Assumption \ref{assumption uniqueness} is satisfied.

\begin{lemma}\label{lemma: unif lip assumption}
    Let $b:[0,T]\times \R^d \times \PP(\R^d) \to \R^d$ and $a:[0,T]\times \R^d \times \PP(\R^d) \to \operatorname{Sym}_+(\R^{d\times d})$ be Borel measurable maps. Assume that $b$ is bounded and that there exists a bounded map $\sigma:[0,T]\times \R^d \times \PP(\R^d) \to \R^{d\times d}$ such that $a = \sigma \sigma^\top$. Moreover, assume that there exists $L>0$ such that
    \begin{equation}
        |\sigma(t,x,\mu) - \sigma(t,y,\mu)| + |b(t,x,\mu) - b(t,y,\mu)| \leq L|x-y| ,
    \end{equation}
    for all $(t,\mu)\in [0,T]\times\PP(\R^d)$ and $x,y\in \R^d$. Then $b$ and $a$ satisfy Assumption \ref{assumption uniqueness}.
\end{lemma}

\begin{proof}
    It is an immediate consequence of \cite[Theorem 6.4]{stroock1997multidimensional}.
\end{proof}

More generally, the assumptions of the previous lemma can be relaxed to whatever case implies that for all $\boldsymbol{\mu}\in C_T(\PP(\R^d))$, the coefficients defined by the maps $(t,x)\mapsto b(t,x,\mu_t)$ and $(t,x)\mapsto a(t,x,\mu_t)$ imply the well-posedness of the Kolmogorov-Fokker-Planck equation associated with them. 
Finally, we can refine Lemma \ref{lemma: transfer of uniqueness}. 

\begin{prop}\label{prop: transfer of uniqueness}
     Let $b:[0,T]\times \R^d \times \PP(\R^d) \to \R^d$ and $a:[0,T]\times \R^d \times \PP(\R^d) \to \operatorname{Sym}_+(\R^{d\times d})$ be Borel measurable maps. Assume they satisfy Assumption \ref{assumption uniqueness}. Then the following are equivalent:
     \begin{enumerate}
        \item for all $s\in[0,T]$ and $\overline{M}\in \PP(\PP(\R^d))$ there exists a unique curve $\boldsymbol{M} = (M_t)_{t\in[s,T]} \in C([s,T],\PP(\PP(\R^d)))$ satisfying $\partial_tM_t = \mathcal{K}_{b_t,a_t}^*M_t$ in $[s,T]$, the integrability condition in \eqref{eq: equation main theorem} (integrating between $s$ and $T$) and $M_s = \overline{M}$;
        \item for all $s\in[0,T]$ and for all $\overline{M} \in \PP(\PP(\R^d))$ there exists a unique $\Lambda \in \PP(C([s,T], \PP(\R^d)))$ concentrated over $\operatorname{KFP}_{[s,T]}(b,a)$, and satisfying \eqref{eq: integr of Lambda} (integrating between $s$ and $T$) and $(\mathfrak{e}_s)_\sharp \Lambda = \overline{M}$;
        \item for all $s\in[0,T]$ and for all $\overline{M} \in \PP(\PP(\R^d))$ there exists a unique $\mathfrak{L} \in \PP(\PP(C([s,T],\R^d)))$ concentrated over $\operatorname{MP}_{[s,T]}(b,a)$, and satisfying \eqref{eq: integr of mathfrak L} (integrating between $s$ and $T$) and $(E_s)_\sharp \mathfrak{L} = \overline{M}$.
    \end{enumerate}
    Moreover, if one of the above holds, then:
    \begin{enumerate}
        \item[(a)] for all $s\in[0,T]$ and $\overline{\mu}\in \PP(\R^d)$ there exists a unique curve of measures $\boldsymbol{\mu} = (\mu_t)_{t\in[s,T]} \in \operatorname{KFP}_{[s,T]}(b,a)$, with $\mu_s = \overline{\mu}$;
        \item[(b)] for all $s\in[0,T]$ and $\overline{\mu}\in \PP(\R^d)$ there exists a unique $\lambda \in \operatorname{MP}_{[s,T]}(b,a)$, with $(\e_s)_\sharp \lambda = \overline{\mu}$.
    \end{enumerate}
\end{prop}

\begin{proof}
    By Lemma \ref{lemma: transfer of uniqueness}, (1)$\iff$(2) and (3)$\implies$(2) are trivial. 

    (2)$\implies$(3): the existence is a consequence of Theorem \ref{thm: from Lambda to mathfrak L}, so we focus on uniqueness. Let $\mathfrak{L}^1,\mathfrak{L}^2 \in \PP(\PP(C([s,T],\R^d)))$ be satisfying the requirements in (3), and define $\Lambda^i := E_\sharp \mathfrak{L}^i$ for $i=1,2$. Note that $E$ is defined as in \eqref{eq: maps for marginals}, but for curves defined in $[s,T]$ instead of $[0,T]$. Thanks to Proposition \ref{prop: from L to Lambda nd}, they both satisfy the requirements in (2), implying that $\Lambda^1 = \Lambda^2$, which will be denoted simply $\Lambda$. Now, consider the disintegrations with respect to $E$
    \[\mathfrak{L}^i = \int_{C([s,T],\PP(\R^d))} \mathfrak{L}_{\boldsymbol{\mu}}^i d\Lambda(\boldsymbol{\mu}), \quad i=1,2,\]
    where $\mathfrak{L}^i_{\boldsymbol{\mu}} \in \PP(\PP(C([s,T],\R^d)))$ is concentrated on $\operatorname{MP}_{[s,T]}(b,a)$ and for $\mathfrak{L}^i$-a.e. $\lambda$, it holds $E(\lambda) = \boldsymbol{\mu}$. We conclude showing there exists a unique element $\lambda$ that satisfies both the properties. Consider $\lambda^1,\lambda^2\in \operatorname{MP}_{[s,T]}(b,a)$, with $E( \lambda^1) = E ( \lambda^2) = \boldsymbol{\mu}$. Since their one-dimensional marginal distributions coincide, they are martingale solutions with respect to the same operator $L_{b_t ,a_t}^{\boldsymbol{\mu}}$, and we denote $\boldsymbol{\mu} = (\mu_t)_{t\in [s,T]} \in \operatorname{KFP}_{[s,T]}(b,a)$. Now, we prove by induction on $n\in \N$ that 
    $\forall s\leq t_1 < \dots < t_n \leq T$, $\forall A_1,\dots,A_n \in \mathcal{B}(\R^d)$
    \begin{equation}\label{eq: induction}
    \lambda^1(\e_{t_1}\in A_1,\dots, \e_{t_n}\in A_n) = \lambda^2(\e_{t_1}\in A_1,\dots, \e_{t_n}\in A_n),
    \end{equation}
    that will make us conclude that $\lambda^1 = \lambda^2$. When $n=1$, \eqref{eq: induction} follows by the fact that the marginals of $\lambda^1$ and $\lambda^2$ coincide. Now, consider $s\leq t_1 < \dots < t_n<t_{n+1} \leq T$ and $\forall A_1,\dots,A_n,A_{n+1} \in \mathcal{B}(\R^d)$. By induction, 
    \[\alpha:= \lambda^1(\e_{t_1}\in A_1,\dots, \e_{t_n}\in A_n) = \lambda^2(\e_{t_1}\in A_1,\dots, \e_{t_n}\in A_n).\]
    If $\alpha=0$, then we are done. Otherwise, define the function $\rho:C([s,T];\R^d) \to [0,+\infty)$ as
    \[\rho(\gamma):= \frac{1}{\alpha}\prod_{i=1}^n \mathds{1}_{A_i}(\e_{t_i}(\gamma)),\]
    that corresponds to the density of $\lambda^1$ and $\lambda^2$ conditioned with respect to $\{\e_{t_1}\in A_1,\dots, \e_{t_n}\in A_n\}$.
    Notice that $\rho$ is $\mathcal{F}_{t_n}$-measurable and consider the measures $\eta^1:=(|_{[t_n,T]})_\sharp (\rho \lambda^1)$ and $\eta^2:= (|_{[t_n,T]})_\sharp (\rho \lambda^2)$. We can conclude that $(\e_{t_n})_\sharp \eta^1 = (\e_{t_n})_\sharp \eta^2$, exploiting the inductive hypothesis \eqref{eq: induction} with a general Borel set $B \in \mathcal{B}(\R^d)$ in place of $A_n$. Moreover, $(\e_{t})_\sharp \eta^1$ and $(\e_{t})_\sharp \eta^2$ are absolutely continuous w.r.t. $\mu_t$ for all $t\in[t_n,T]$ and their densities are controlled by $1/\alpha$. Then, thanks to Lemma \ref{lemma: restriction} and Assumption \ref{assumption uniqueness}, we can assess that $(\e_{t})_\sharp (\rho \lambda^1) = (\e_{t})_\sharp (\rho \lambda^2)$ for all $t\in [t_n,T]$. In particular, it holds
    \begin{align*}
        & \hspace{-1cm} \frac{\lambda^1(\e_{t_1}\in A_1,\dots,\e_{t_n}\in A_n, \e_{t_{n+1}}\in A_{n+1})}{\alpha} = \big[(\e_{t_{n+1}})_\sharp (\rho \lambda^1)\big](A_{n+1})
        \\
        =&  \big[(\e_{t_{n+1}})_\sharp (\rho \lambda^2)\big](A_{n+1}) = \frac{\lambda^2(\e_{t_1}\in A_1,\dots,\e_{t_n}\in A_n, \e_{t_{n+1}}\in A_{n+1})}{\alpha}.
    \end{align*}
    \noindent As in Lemma \ref{lemma: transfer of uniqueness}, conditions (a) and (b) are byproduct of the previous arguments, or simply follows considering $\overline{M} = \delta_{\overline{\mu}}$ as initial condition. 
\end{proof}

\normalcolor

\begin{oss}
    The previous proposition can be localized, in the sense that we may substitute $\overline{M} \in \PP(\PP(\R^d))$ and $\overline{\mu} \in \PP(\R^d)$, respectively, with $\overline{M}$ concentrated over a given set $\mathcal{A}\subset \PP(\R^d)$ and $\overline{\mu} \in \mathcal{A}$.
\end{oss}

\bigskip

\noindent \textbf{Acknowledgments.} The author warmly thanks Giuseppe Savar\'e and Dario Trevisan for their fruitful suggestions and the encouragement in pursuing the results presented in this paper.

\printbibliography

\end{document}